\documentclass[11pt, reqno]{amsart}
\usepackage{amsfonts,latexsym}
\usepackage{amsmath}
\usepackage{amscd}
\usepackage{float,amsmath,amssymb,mathrsfs,bm,multirow,graphics}
\usepackage[dvips]{graphicx}
\usepackage[percent]{overpic}
\usepackage[numbers,sort&compress]{natbib}
\usepackage{enumerate}

\newcommand{\nequation}{\setcounter{equation}{0}}

\newcommand{\R}{{\Bbb R}}

\newcommand{\C}{{\Bbb C}}

\newcommand{\proofbegin}{\noindent{\it Proof.\,\,}}
\newcommand{\proofend}{\hfill$\Box$\bigskip}

\DeclareMathOperator{\im}{Im}
\DeclareMathOperator{\re}{Re}

\def\XXint#1#2#3{{\setbox0=\hbox{$#1{#2#3}{\int}$}
\vcenter{\hbox{$#2#3$}}\kern-.5\wd0}}

\newtheorem{theorem}{Theorem}[section]

\newtheorem{lemma}[theorem]{Lemma}
\newtheorem{definition}[theorem]{Definition}

\newtheorem{remark}[theorem]{Remark}

\newtheorem{figuretext}{Figure}

\input epsf
\title[Admissible boundary values]{\sc Admissible boundary values for the defocusing nonlinear Schr\"odinger equation with asymptotically time-periodic data}

\author{Jonatan Lenells}
\address{Department of Mathematics, KTH Royal Institute of Technology, 100 44 Stockholm, Sweden.}
\email{jlenells@kth.se}

\begin{document}

\begin{abstract} 
\noindent
We consider solutions of the defocusing nonlinear Schr\"odinger equation in the quarter plane whose Dirichlet boundary data approach a single exponential $\alpha e^{i\omega t}$ as $t \to \infty$. In order to determine the long time asymptotics of the solution, it is necessary to first characterize the asymptotic behavior of the Neumann value in terms of the given data. Assuming that the initial data decay as $x \to \infty$, we derive necessary conditions for the Neumann value to asymptote towards a single exponential of the form $ce^{i\omega t}$. Since our approach yields expressions which relate $\alpha$, $\omega$, and $c$, the result can be viewed as a characterization of the large $t$ behavior of the Dirichlet to Neumann map for single exponential profiles.
\end{abstract}

\maketitle

\noindent
{\small{\sc AMS Subject Classification (2000)}: 37K15, 35Q15.}

\noindent
{\small{\sc Keywords}: Initial-boundary value problem, integrable system, long-time asymptotics.}

%\tableofcontents

\section{Introduction}\nequation
The long time asymptotics of solutions of integrable PDEs can be analyzed by means of the nonlinear steepest descent method for Riemann-Hilbert (RH) problems introduced by Deift and Zhou \cite{DZ1993} in the 1990s. In the context of initial value problems, the RH problem is formulated in terms of certain spectral functions whose definitions involve the initial data of the solution. In this way, the asymptotics of the solution of the initial value problem for the modified KdV \cite{DZ1993}, the nonlinear Schr\"odinger (NLS) \cite{DIZ1993}, and other integrable equations \cite{BS2008, DVZ1994, DZ1995, GT2009, K1993} has been rigorously established. 

In the case of initial-boundary value (IBV) problems, the picture is much less complete. In this case, the RH problem is formulated in terms of spectral functions whose definitions involve the initial data as well as the boundary values of the solution. Under the assumption that both the Dirichlet and Neumann boundary values decay as $t \to \infty$, the long time asymptotics of the solution on the half-line has been determined for several equations, such as the modified KdV \cite{BFS2004} and NLS \cite{FI1996, FIS2005} equations, using the nonlinear steepest descent method together with the unified transform formalism introduced in \cite{F1997}. 
If the boundary data do not decay as $t \to\infty$, the situation is significantly more complicated due to the fact that only a subset of the boundary values are known for a well-posed problem. For example, for the Dirichlet problem for the NLS equation on the half-line, the initial and Dirichlet data are known, whereas the Neumann data have to be determined as part of the solution. 
Thus, in order to successfully apply the Deift-Zhou approach to this problem, it appears necessary to first determine the long time asymptotics of the Neumann value. Even though it is possible to characterize the Neumann value in terms of the initial and Dirichlet data via a system of nonlinear integral equations \cite{BFS2003, F2005, trilogy1}, this is not sufficient for finding the long time asymptotics. 

From the point of view of applications, an important class of IBV problems involving non-decaying boundary data consists of problems with a time-periodic (or at least asymptotically time-periodic) boundary condition. This type of problem arises, for example, for the KdV equation when modeling the near-shore wave motion generated by waves propagating from deep water, or when studying waves in a wave tank with a periodically moving wave maker mounted at one end \cite{BPS1981}.
In the case of asymptotically time-periodic boundary data, the construction of the Dirichlet to Neumann map simplifies and it is in fact possible to find the asymptotic form of the unknown Neumann value for the Dirichlet problem of the NLS and other integrable equations perturbatively to all orders by means of a recursive scheme \cite{FLtperiodic2}. 
If we restrict the class of boundary conditions even further and consider the special case of a Dirichlet value which asymptotes towards a {\it single} exponential,
\begin{align}\label{singleexp}
q(0,t) \sim \alpha e^{i\omega t}, \qquad t \to \infty, \quad \alpha > 0, \quad \omega \in \R,
\end{align}
then pioneering asymptotic formulas have been established for the focusing NLS equation in a series of papers by Boutet de Monvel and coauthors \cite{BIK2007, BIK2009, BKS2009, BKSZ2010}. Defining a pair of functions $\{g_0^b(t), g_1^b(t)\}$ as {\it asymptotically admissible} if there exists a solution $q(x,t)$ of the NLS in the quarter plane $\{x>0, t>0\}$ such that $q$ decays as $x \to \infty$ and
$$q(0, t) \sim g_0^b(t), \qquad q_x(0, t) \sim g_1^b(t), \qquad t \to \infty,$$
they show that the pair $\{\alpha e^{i\omega t}, ce^{i\omega t}\}$ where $\alpha > 0$, $\omega \in \R$, and $c \in \C$, is asymptotically admissible if and only if the parameters $(\alpha,\omega,c)$ satisfy either
\begin{subequations}\label{1.3}
\begin{align}\label{1.3a}
 c = \pm \alpha\sqrt{\omega - \alpha^2} \quad \text{and} \quad \omega \geq \alpha^2,
\end{align}
or
\begin{align}\label{1.3b}
  c = i \alpha \sqrt{|\omega| + 2 \alpha^2} \quad \text{and}  \quad \omega \leq - 6\alpha^2.
\end{align}
\end{subequations}
Furthermore, if one of the conditions (\ref{1.3}) is fulfilled, they determine the long-time asymptotics of $q(x,t)$ by applying the nonlinear steepest descent method.

Our goal in this paper is to take a first few steps towards deriving the analogs of the results of \cite{BIK2007, BIK2009, BKS2009, BKSZ2010} in the case of the {\it defocusing} NLS. Even for the relatively simple example of a single exponential, the analysis of the defocusing NLS is surprisingly rich. Indeed, in the spectral space, the admissibility of a given set of initial and boundary values is encoded in the so-called global relation, which imposes certain conditions on the spectral functions. For the focusing NLS, the imposition of an appropriate analyticity condition derived from the global relation leads to the two families of admissible parameter triples in (\ref{1.3}). 
On the other hand, for the defocusing NLS, imposing the analogous analyticity condition leads to {\it five} different families of triples; in addition to the analogs of the branches (\ref{1.3}) present in the focusing case, there are three branches for which both the real and the imaginary parts of $c$ are nonzero. Since each of these additional branches depends on two or three parameters, this provides a large number of potentially asymptotically admissible Dirichlet and Neumann pairs for the defocusing NLS.

Our main result (Theorem \ref{mainth}) is presented in Section \ref{mainsec}. After having introduced appropriate eigenfunctions and spectral functions in Section \ref{eigensec}, the proof is given in Section \ref{proofsec}. In Section \ref{analysissec}, we analyze the five families found in the main theorem further. 

\section{Main Result}\nequation\label{mainsec}
We consider the NLS equation
\begin{align}\label{nls}
  iq_t + q_{xx} - 2\lambda |q|^2 q = 0, \qquad \lambda = \pm 1,
\end{align}
in the quarter plane $\{x> 0, t> 0\}$. The two versions of (\ref{nls}) with $\lambda = -1$ and $\lambda = 1$ are referred to as the focusing and defocusing NLS respectively. We let $\mathcal{S}([0,\infty))$ denote the Schwartz class
\begin{align}\label{schwartzdef}
\mathcal{S}([0,\infty)) = \{u \in C^\infty([0,\infty)) \, | \, x^n u^{(m)}(x) \in L^\infty([0,\infty)) \text{ for all } n, m \geq 0\}.
\end{align}

\begin{definition}\upshape\label{soldef}
A {\it solution of the NLS in the quarter plane} is a smooth function $q:[0,\infty) \times [0,\infty) \to \C$ such that $q(\cdot, t) \in \mathcal{S}([0,\infty))$ for each $t \geq 0$, and such that (\ref{nls}) is satisfied for $x> 0$ and $t > 0$.
\end{definition}

\begin{definition}\upshape\label{admissibledef}
A pair of functions $\{g_0^b(t), g_1^b(t)\}$, defined for $t \geq 0$, is {\it asymptotically admissible} for NLS if there exists a solution $q(x,t)$ of the NLS in the quarter plane such that the Dirichlet and Neumann boundary values of $q$ asymptote towards $g_0^b(t)$ and $g_1^b(t)$ respectively in the sense that
\begin{align*}
q(0, t) - g_0^b(t) = O(t^{-7/2}), \qquad q_x(0,t) - g_1^b(t) = O(t^{-7/2}), \qquad t\to \infty.
\end{align*}
\end{definition}

In the context of the single exponential profile (\ref{singleexp}), it is convenient to introduce the notion of an admissible parameter triple as follows.

\begin{definition}\upshape
The parameter triple $(\alpha, \omega, c)$ where $\alpha>0$, $\omega \in \R$, and $c \in \C$ is {\it admissible} if the pair $\{\alpha e^{i\omega t}, ce^{i\omega t}\}$ is asymptotically admissible.
\end{definition}

It was shown in \cite{BKS2009} that the admissible triples for the focusing NLS are given by (\ref{1.3}). 
The following theorem determines potentially admissible triples in the defocusing case.

\begin{theorem}\label{mainth}
Every admissible triple $(\alpha, \omega, c)$ for the defocusing NLS equation where $\alpha>0$, $\omega \in \R$, and $c \in \C$, belongs to one of the following disjoint subsets:
\begin{subequations}\label{admissiblesets}
\begin{align}\label{admissibleA}
& \bigg\{\bigg(\alpha, \omega, c = \pm \sqrt{\frac{(\omega + 3\alpha^2)^3}{27\alpha^2}} + \frac{i|\omega|^{3/2}}{3\sqrt{3} \alpha}\bigg) \; \bigg| \;  \alpha > 0, \; -3\alpha^2 \leq \omega < 0\bigg\},
	\\ \nonumber
& \bigg\{\bigg(\alpha  = -\frac{4K^3 + \omega K}{c_2}, \omega, c = \pm \sqrt{ \bigg(\alpha^2+ \frac{\omega}{2}\bigg)^2 -c_2^2 - 2K^2(6K^2 + \omega)} + ic_2\bigg)  
	\\ \label{admissibleB}
&\hspace{3.5cm} \bigg| \;  -12 K^2 < \omega < -4K^2, \; 0 < c_2 \leq -\frac{4K^2 + \omega}{2}, \; K > 0\bigg\}, 
	\\ \label{admissibleC}
& \big\{(\alpha,\omega, c = i\alpha\sqrt{-2\alpha^2 - \omega}) \; \big| \; \alpha > 0, \; \omega < -3\alpha^2\big\},
	\\ \label{admissibleD}
& \big\{(\alpha,\omega,c = \pm \alpha\sqrt{\omega + \alpha^2}) \; \big| \; \omega + \alpha^2 \geq 0, \; \alpha > 0\big\},
	\\ \nonumber
& \bigg\{\bigg(\alpha = -\frac{4K^3 + \omega K}{c_2}, \omega, c = \pm \sqrt{\bigg(\alpha^2+ \frac{\omega}{2}\bigg)^2 -c_2^2 - 2K^2(6K^2 + \omega)} + ic_2\bigg)
	\\ \label{admissibleE}
& \hspace{3.5cm} \bigg| \; -4K^2 < \omega \leq -3K^2, \; -\frac{4K^2 + \omega}{2} \leq c_2 < 0, \; K > 0\bigg\}.
\end{align}
\end{subequations}
\end{theorem}

\begin{remark}\upshape
The subset of (\ref{admissibleA}) for which $\omega = -3\alpha^2$ can be combined with (\ref{admissibleC}) to give the following set of triples:
\begin{align}\label{tripleplane}
\big\{(\alpha,\omega, c = i\alpha\sqrt{-2\alpha^2 - \omega}) \; \big| \; \alpha > 0, \; \omega \leq -3\alpha^2\big\}.
\end{align}
This set is the defocusing analog of the branch of admissible triples (\ref{1.3b}) present in the focusing case. Similarly, (\ref{admissibleD}) is the analog of (\ref{1.3a}). The branches (\ref{admissibleA}), (\ref{admissibleB}), and (\ref{admissibleE}) are new to the defocusing case and have no focusing analogs.
\end{remark}

%In case (\ref{case32}), we have $c_1 = 0$ iff $K = \alpha/2$. The two-parameter subfamily for which $K = \alpha/2$ can be written as
%$$\{(\alpha,\omega, ic_2) | \alpha > 0, \omega = -\alpha^2 - 2 c_2, c_2 \in (0, \alpha^2)\}.$$

\begin{remark}\upshape
  Theorem \ref{mainth} states that the parameter set $(\alpha, \omega, c)$ of any asymptotically admissible pair of the form $\{\alpha e^{i\omega t}, ce^{i\omega t})$ must satisfy one of the conditions in (\ref{admissiblesets}). However, it remains an open problem to determine whether all triples in (\ref{admissiblesets}) are actually admissible. Although the condition we use to find (\ref{admissiblesets}) below is a necessary condition for admissibility, it is not guaranteed that it is also sufficient.
In the case of single exponential profiles for the focusing NLS, sufficiency was established by using an associated RH problem to construct a solution whose Dirichlet and Neumann values have the desired asymptotics \cite{BK2007, BKS2009}. 
  \end{remark}

\section{Eigenfunctions and spectral functions}\label{eigensec}\nequation
Equation (\ref{nls}) admits the Lax pair
\begin{align}\label{lax}
\begin{cases}
  \phi_x + ik\sigma_3 \phi = Q\phi,
  	\\ 
  \phi_t + 2ik^2 \sigma_3 \phi = \tilde{Q} \phi,	
\end{cases}
\end{align}
where $k \in \C$ is the spectral parameter, $\phi(x,t,k)$ is a $2\times 2$-matrix valued eigenfunction, and
\begin{align}\label{QQtildedef}
Q = \begin{pmatrix} 0 & q \\
\lambda \bar{q} & 0 \end{pmatrix}, \qquad
\tilde{Q} = \begin{pmatrix} -i\lambda |q|^2 & 2kq + iq_x \\
2k \lambda \bar{q} - i\lambda \bar{q}_x & i\lambda |q|^2 \end{pmatrix}, \qquad \sigma_3 = \begin{pmatrix} 1 & 0 \\ 0 & -1 \end{pmatrix}.
\end{align}

Suppose $(\alpha, \omega, c)$ is an admissible triple and let $q(x,t)$ be an associated solution of the NLS in the quarter plane such that
\begin{align*}
  q(0,t) - \alpha e^{i\omega t} = O(t^{-7/2}), \qquad q_x(0, t) - ce^{i\omega t} = O(t^{-7/2}), \qquad t\to \infty.
\end{align*}

\subsection{The background eigenfunction}
We define the `background' eigenfunction $\psi^b(t,k)$ by
\begin{align}\label{backgroundefunction}
\psi^b(t,k) = e^{\frac{i\omega}{2} t \sigma_3} E(k) e^{-i\Omega(k)t\sigma_3},
\end{align}
where $\Omega(k)$ and $E(k)$ are defined by
\begin{align} \label{Omega2def}
& \Omega(k) = \sqrt{4k^4 + 2\omega k^2 + 4\lambda \alpha \im(c) k + \bigg(\frac{\omega}{2} + \lambda \alpha^2\bigg)^2 - \lambda |c|^2},
	\\ \label{Edef}
& E(k) = \sqrt{\frac{2\Omega - H}{2\Omega}} \begin{pmatrix} 1 & \frac{\lambda iH}{2\alpha k - i\bar{c}} \\
-\frac{iH}{2\alpha k + ic} & 1 \end{pmatrix},	
\end{align}
with 
$$H(k) = \Omega(k) - 2k^2 - \lambda \alpha^2 - \frac{\omega}{2}.$$
The function $\psi^b(t,k)$ is a solution of the background $t$-part 
\begin{align}\label{backgroundtpart}
  \psi_t^b + 2ik^2\sigma_3 \psi^b = \tilde{Q}^b \psi^b
\end{align}
where $\tilde{Q}^b$ is given by the expression (\ref{QQtildedef}) for $\tilde{Q}$ with $q$ and $q_x$ replaced with $\alpha e^{i\omega t}$ and $ce^{i\omega t}$ respectively.
We view the above functions as being defined on the cut complex $k$-plane $\C \setminus \mathcal{C}$, where $\mathcal{C}$ is the union of two branch cuts connecting the four roots of $\Omega^2(k)$ and a third branch cut which connects the points $\frac{i\bar{c}}{2\alpha}$ and $-\frac{ic}{2\alpha}$. We note in this regard that the equality
\begin{align}\label{2OmegaHH}
  (H - 2\Omega)H = \lambda (2\alpha k - i\bar{c})(2\alpha k + ic)
\end{align}  
implies that the zeros of $2\Omega - H$ and $H$ are included in the set $\{\frac{i\bar{c}}{2\alpha}, -\frac{ic}{2\alpha}\}$.
The branches of $\Omega(k)$ and $\sqrt{\frac{2\Omega -H}{2\Omega}}$ are fixed by their asymptotics as $k \to \infty$:
\begin{align}\label{Omegasqrt}
\Omega(k) = 2k^2 + \frac{\omega}{2} + O(k^{-1}), \qquad
\sqrt{\frac{2\Omega - H}{2\Omega}} = 1 + O(k^{-1}), \qquad k \to \infty.
\end{align}
Note that $\det E(k) = 1$. We will assume that the branch cuts in $\mathcal{C}$ intersect transversely in at most finitely many points. 
The function $\Omega(k)$ changes sign as $k$ crosses one of the branch cuts connecting two of the zeros of $\Omega^2$, whereas it does not jump across the branch cut connecting $\frac{i\bar{c}}{2\alpha}$ and $-\frac{ic}{2\alpha}$.

\subsection{The eigenfunctions $\{\phi_j\}_1^3$}
Let $\hat{\sigma}_3$ act on a $2\times 2$ matrix $A$ by $\hat{\sigma}_3A = [\sigma_3, A]$, so that $e^{\hat{\sigma}_3} A = e^{\sigma_3} A e^{-\sigma_3}$.
We introduce three solutions $\{\phi_j(x,t,k)\}_1^3$ of (\ref{lax}) by
\begin{align}\nonumber
 & \phi_1(x,t,k) = \mu_1(x,t,k) e^{-i(kx + (\Omega(k) - \frac{\omega}{2})t)\sigma_3},
  	\\ \label{phidef}
 &  \phi_j(x,t,k) = \mu_j(x,t,k) e^{-i(kx + 2k^2 t)\sigma_3}, \qquad j = 2,3,
\end{align}
where the $\mu_j$'s are $2\times 2$-matrix valued solutions of the Volterra integral equations
\begin{align}\nonumber
\mu_1(x,t,k) = &\; e^{-ikx\hat{\sigma}_3}\bigg\{  \mathcal{E}(t,k) 
- \mathcal{E}(t,k) \int_t^\infty e^{i(\Omega(k) - \frac{\omega}{2})(t' - t)\hat{\sigma}_3} \big[\mathcal{E}^{-1}(t', k) 
	\\\nonumber
&\times (\tilde{Q} - \tilde{Q}^b)(0,t',k) \mu_1(0,t',k) \big] dt' 
+ \int_0^x e^{ikx' \hat{\sigma}_3}[Q(x',t) \mu_1(x',t,k)] dx'\bigg\}
	\\ \label{mudef}
\mu_j(x,t,k) = &\; I + \int_{(x_j, t_j)}^{(x,t)} e^{i[k(x'-x) + 2k^2(t'-t)] \hat{\sigma}_3} W_j(x',t',k), \qquad j = 2,3,
\end{align}
with $(x_2, t_2) = (0,0)$, $(x_3, t_3) = (\infty, t)$, and
$$\mathcal{E}(t,k) = e^{\frac{i\omega}{2}t\hat{\sigma}_3}E(k), \qquad 
W_j = (Q dx + \tilde{Q} dt)\mu_j, \quad j =2,3.$$

Let
\begin{align*}
D_1 = \{\im k > 0\} \cap \{\im \Omega(k) > 0\},  \qquad
D_2 = \{\im k > 0\} \cap \{\im \Omega(k) < 0\}, 
	\\
D_3 = \{\im k < 0\} \cap \{\im \Omega(k) > 0\},  \qquad
D_4 = \{\im k < 0\} \cap \{\im \Omega(k) < 0\},
\end{align*}
and let $D_+ = D_1 \cup D_3$ and $D_- = D_2 \cup D_4$.
The asymptotics (\ref{Omegasqrt}) of $\Omega(k)$ implies that $D_j$ can be viewed as a deformation of the $j$'th quadrant of the complex $k$-plane. 
The eigenfunctions $\{\mu_j(x,t,k)\}_1^3$ have the following properties (see \cite{FLtperiodic} for a proof in a more general context).

\begin{itemize}
  \item The first (resp. second) column of $\mu_1(0,t,k)$ is defined and analytic for $D_- \setminus \mathcal{C}$ (resp. $D_+ \setminus \mathcal{C}$).
The second column of $\mu_1$ has a continuous extension to the boundary of $D_+ \setminus \mathcal{C}$ in the sense that away from the branch points, the limits from the right and left onto every branch cut in $D_+$, and the limits onto each part of the boundary of $D_+$, exist and are continuous. If a branch cut can be approached from both right and left from within $D_+ \setminus \mathcal{C}$, then the right and left limits are, in general, different.
%Near the branch points where $\tr Z = \pm 2$ (i.e. where $z(k) = \pm 1$ and hence $\tilde{\Omega}(k) = \pi n/T$, $n \in \Z$), $\mu_1(x,t,k)$ generically has the same type of singularity as $\mathcal{E}(t,k)$ (generically of the form $(\tr Z \mp 2)^{-1/4}$).
   
  \item $\mu_2(x,t,k)$ is defined and analytic for all $k \in \C$.

  \item The first (resp. second) column of $\mu_3(x,t,k)$ is defined and analytic for $\im k < 0$ (resp. $\im k > 0$) with a continuous extension to $\im k \leq 0$ (resp. $\im k \geq 0$).

\item The $\mu_j$'s are normalized so that
\begin{align*}
  \lim_{t \to \infty}[\mu_1(0,t,k) - \mathcal{E}(t,k)] & = 0, \qquad k \in (D_-\setminus \mathcal{C}, D_+\setminus \mathcal{C}),
  	\\
 \mu_2(0,0,k)&  = I, \qquad k \in \C,
	\\
  \lim_{x \to \infty} \mu_3(x,0,k)&  = I, \qquad k \in (\bar{\C}_-, \bar{\C}_+),
\end{align*}
where the notation $k \in (A_1, A_2)$ indicates that the first and second columns are valid for $k \in A_1$ and $k \in A_2$, respectively. More precisely, if $K_\pm$ are compact subsets of $(\overline{D_\pm \setminus \mathcal{C}})\setminus \mathcal{P}$, where $\mathcal{P}$ denotes the set of branch points, then
\begin{align}\label{mu1calE}
 |\mu_1(0,t,k) - \mathcal{E}(t,k)| \leq C(1+t)^{-5/2}, \qquad k \in (K_-, K_+), \quad t \geq 0.
\end{align}
\end{itemize}

\subsection{Spectral functions}
We define the spectral functions $s(k)$ and $S(k)$ by
$$s(k) = \mu_3(0,0,k), \qquad S(k) = \mu_1(0,0,k),$$
and write
$$s(k) = \begin{pmatrix} \overline{a(\bar{k})} & b(k) \\
\lambda \overline{b(\bar{k})} & a(k) \end{pmatrix}, \qquad
S(k) = \begin{pmatrix} \overline{A(\bar{k})} & B(k) \\
\lambda \overline{B(\bar{k})} & A(k) \end{pmatrix}.$$
Then
$$\phi_3(x,t,k) = \phi_2(x,t,k) s(k), \qquad
\phi_1(x,t,k) = \phi_2(x,t,k) S(k).$$

The analyticity properties of $\mu_1$ and $\mu_3$ imply corresponding analyticity properties for the spectral functions. In particular, $A(k)$ and $B(k)$ are defined and analytic for $k \in D_+ \setminus \mathcal{C}$ with a continuous extension to $\bar{D}_+ \setminus \mathcal{C}$. Moreover, away from the branch points, $A(k)$ and $B(k)$ have continuous extensions onto any branch cut that intersects $\bar{D}_+$. If $k\in D_+$ can approach the branch cut from both the left and right sides, we denote the corresponding limits, which in general are different, by $\{A_-,B_-\}$ and $\{A_+,B_+\}$ respectively.
%Furthermore, $A(k) \to 1$ and $B(k) \to 0$ as $k \to \infty$, $k \in \bar{D}_+ \setminus \mathcal{C}$.
The functions $a(k)$ and $b(k)$  are defined and analytic in $\im k > 0$ with a continuous extension to $\im k \geq 0$.

\subsection{The global relation}
Letting $T \to \infty$ in the $(12)$ entry of the equation
$$S^{-1}(k)s(k) = e^{i(\Omega(k) - \frac{\omega}{2})T \sigma_3}(\mu_1^{-1}(0,T,k)\mu_3(0,T,k))e^{-2ik^2T\sigma_3}$$
and using the decay of the exponential $e^{i(\Omega(k) + 2k^2)T}$, we find
$$A(k)b(k) - a(k) B(k) = 0, \qquad k \in D_1 \setminus \mathcal{C}, \quad \im (\Omega(k) + 2k^2) > 0.$$
Assuming that $D_1 \setminus \mathcal{C}$ is connected, the condition $\im (\Omega(k) + 2k^2) > 0$ can be removed by analytic continuation. This yields the following global relation:
\begin{align}\label{GR}
A(k) b(k) - a(k) B(k) = 0, \qquad k \in \overline{D_1 \setminus \mathcal{C}} \quad \text{(if $D_1 \setminus \mathcal{C}$ is connected)}.
\end{align}

\subsection{Inadmissible triples}
Using the global relation, we will now show that the jumps of $E(k)$ across the branch cuts connecting the zeros of  $\Omega^2(k)$ lead to jumps also in the quotient $\frac{B(k)}{A(k)}$ along these cuts (see \cite{BKS2009} for the analogous result in the focusing case). 
The following lemma is essential for the proof of Theorem \ref{mainth}.

\begin{lemma}\label{inadmissiblelemma}
Assume that $D_1 \setminus \mathcal{C}$ is connected and that there exists an open set $U \subset \bar{D}_1$ such that one of the two branch cuts connecting the four zeros of $\Omega^2(k)$ intersects $U$. Then the triple $(\alpha, \omega, c)$ is inadmissible. 
\end{lemma}
\proofbegin
Let $C$ be a branch cut connecting two zeros of $\Omega^2(k)$ which intersects $U$. 
Let  $\mathcal{E}_-$ and $\mathcal{E}_+$ denote the limits of $\mathcal{E}$ onto $C \cap U$ from the left and right, respectively. Let $(\mu_1(0,t,k))_{\pm}$ be defined in terms of $\mathcal{E}_\pm(t,k)$ via (\ref{mudef}).
Since $\Omega_+ = - \Omega_-$ on $C$ and $U \subset \bar{D}_1$, we have $\im \Omega = 0$ on $U\cap C$.  Hence both columns of $(\mu_1)_\pm$ are well-defined. The eigenfunctions $\nu_\pm(t,k)$ introduced by
\begin{align}\label{mu1nuE}
  (\mu_1(0,t,k))_\pm = \nu_\pm(t,k) \mathcal{E}_\pm(t,k), \qquad k \in C \cap U,
\end{align}
satisfy the integral equation
\begin{align}\nonumber
\nu_\pm(t,k) = &\; I - \int_t^\infty \psi^b(t,k) (\psi^b)^{-1}(t',k) (\tilde{Q} - \tilde{Q}^b)(0,t', k) \nu_\pm(t', k)
	\\\label{nuVolterra}
&\times \psi^b(t', k) (\psi^b)^{-1}(t, k) dt',
\end{align}
where we note that, since $\tilde{Q}^b$ is a polynomial in $k$, $\psi^b(t,k)(\psi^b)^{-1}(t',k)$ and its inverse are entire functions of $k$. The assumption that $\tilde{Q} - \tilde{Q}^b = O(t^{-5/2})$ implies that the Volterra equation (\ref{nuVolterra}) has a unique solution for $k \in C \cap U$. Hence $\nu_- = \nu_+$.

Evaluating the second column of (\ref{mu1nuE}) at $t = 0$, we find
\begin{align}\label{BApm}
\begin{pmatrix} B(k) \\ A(k) \end{pmatrix}_{\pm} = \nu(0,k) \begin{pmatrix} E_{12}(k) \\ E_{22}(k) \end{pmatrix}_{\pm}.
\end{align}
Using the short-hand notation $\nu_{ij}$ for the $(ij)$'th entry of $\nu(0,k)$, equations (\ref{Edef}) and (\ref{BApm}) imply
\begin{align}\label{fracBA}
\bigg(\frac{B(k)}{A(k)}\bigg)_\pm = \frac{\nu_{11} \lambda i H_{\pm} + \nu_{12}(2\alpha k - i\bar{c})}{\nu_{21} \lambda i H_\pm + \nu_{22}(2\alpha k - i\bar{c})}.
\end{align}
Since $H_+ - H_- = \Omega_+ - \Omega_- = 2\Omega_+$ and $\det \nu = 1$, we find
\begin{align*}
\bigg(\frac{B(k)}{A(k)}\bigg)_+ - \bigg(\frac{B(k)}{A(k)}\bigg)_-
= \frac{- 2 \lambda (2i\alpha k + \bar{c}) \Omega_+ }{(\lambda \nu_{21}H_+ - \nu_{22}(2i\alpha k + \bar{c}))(\lambda \nu_{21} H_- - \nu_{22}(2i\alpha k + \bar{c})} \neq 0.
\end{align*}
It follows that the quotient $B(k)/A(k)$ is discontinuous across $C \cap U$. Since $a(k)$ and $b(k)$ are continuous in $\bar{D}_1 \cup \bar{D}_2$, this contradicts the global relation (\ref{GR}), showing that the triple $(\alpha, \omega, c)$ cannot be admissible.
\proofend

%Let $\mathcal{B} \subset \mathcal{C}$ denote the set of branch points:
%$$\mathcal{B} = \big\{k \in \C | \Omega^2(k) = 0\big\} \cup \bigg\{\frac{i\bar{c}}{2\alpha}, -\frac{ic}{2\alpha}\bigg\}.$$

\section{Proof of Theorem \ref{mainth}}\label{proofsec}\nequation
The proof is inspired by the analogous proof in the focusing case \cite{BKS2009}, but is complicated by the positive sign of $\lambda$.

Suppose $(\alpha, \omega, c)$ is an admissible triple and let $q(x,t)$ be the associated solution of (\ref{nls}) with the properties listed in Definition \ref{admissibledef}. 
We will see that this leads to a contradiction unless $(\alpha, \omega, c)$ belongs to one of the sets in (\ref{admissiblesets}).
Since we are considering the defocusing case, we henceforth set $\lambda = 1$.

Letting $k = k_1 + ik_2$ and $c = c_1 + ic_2$, we can write 
$$\im \Omega^2(k) = 4k_2(4k_1(k_1^2 - k_2^2) + \omega k_1 + \alpha c_2).$$
Hence
$$\{\im \Omega(k) = 0\} = \{\im \Omega^2(k) = 0\} \cap \{\re \Omega^2(k) \geq 0\},$$
where
$$\{\im \Omega^2(k) = 0\}
= \R \cup \Gamma \quad \text{with} \quad \Gamma := \{k \in \C \,|\, 4 k_1 (k_1^2 - k_2^2) + \omega k_1 + \alpha c_2 = 0\}$$
and
\begin{align*}
 \{\re \Omega^2(k) \geq 0\} = \bigg\{k \in \C \,\bigg|\, &4(k_1^2 - k_2^2)^2 - 16 k_1^2 k_2^2 + 2\omega (k_1^2 - k_2^2) + 4 \alpha c_2k_1
	\\
& + \bigg(\frac{\omega}{2} + \alpha^2\bigg)^2 - |c|^2 \geq 0\bigg\}.
\end{align*}
Each of the sets $\{\im \Omega^2(k) = 0\}$ and $\{\re \Omega^2(k) = 0\}$ consists of four curves that asymptote towards the rays $\arg k = \frac{n \pi}{4}$ and $\arg k = \frac{\pi}{8} + \frac{n \pi}{4}$, $n = 0, \dots, 7$, respectively as $|k| \to \infty$. The zeros of $\Omega^2(k)$ are the intersections points of these sets.

\begin{lemma}
$\Omega^2(k)$ has a zero of order $\geq 2$ at a point $K \in \C$ if and only if
\begin{align}\label{doubleconditions}  
4K^3  + \omega K + \alpha c_2= 0 \quad \text{and} \quad
c_1^2 + c_2^2 + 2K^2(6K^2 + \omega) = \bigg(\alpha^2 + \frac{\omega}{2}\bigg)^2.
\end{align}
\end{lemma}
\proofbegin
A fourth order polynomial $P(k) = k^4 + a_2 k^2 + a_1 k + a_0$
has a zero of order $\geq 2$ at $K \in \C$ if and only if $P(K) = P'(K) = 0$, that is, if and only if $K$ satisfies the equations
$$4K^3 + 2a_2 K + a_1 = 0, \qquad K^2(a_2 + 3K^2) = a_0.$$
Letting $P(k) = \Omega^2(k)/4$, we find (\ref{doubleconditions}). 
\proofend 

\begin{figure}
\begin{center}
\begin{overpic}[width=.3\textwidth]{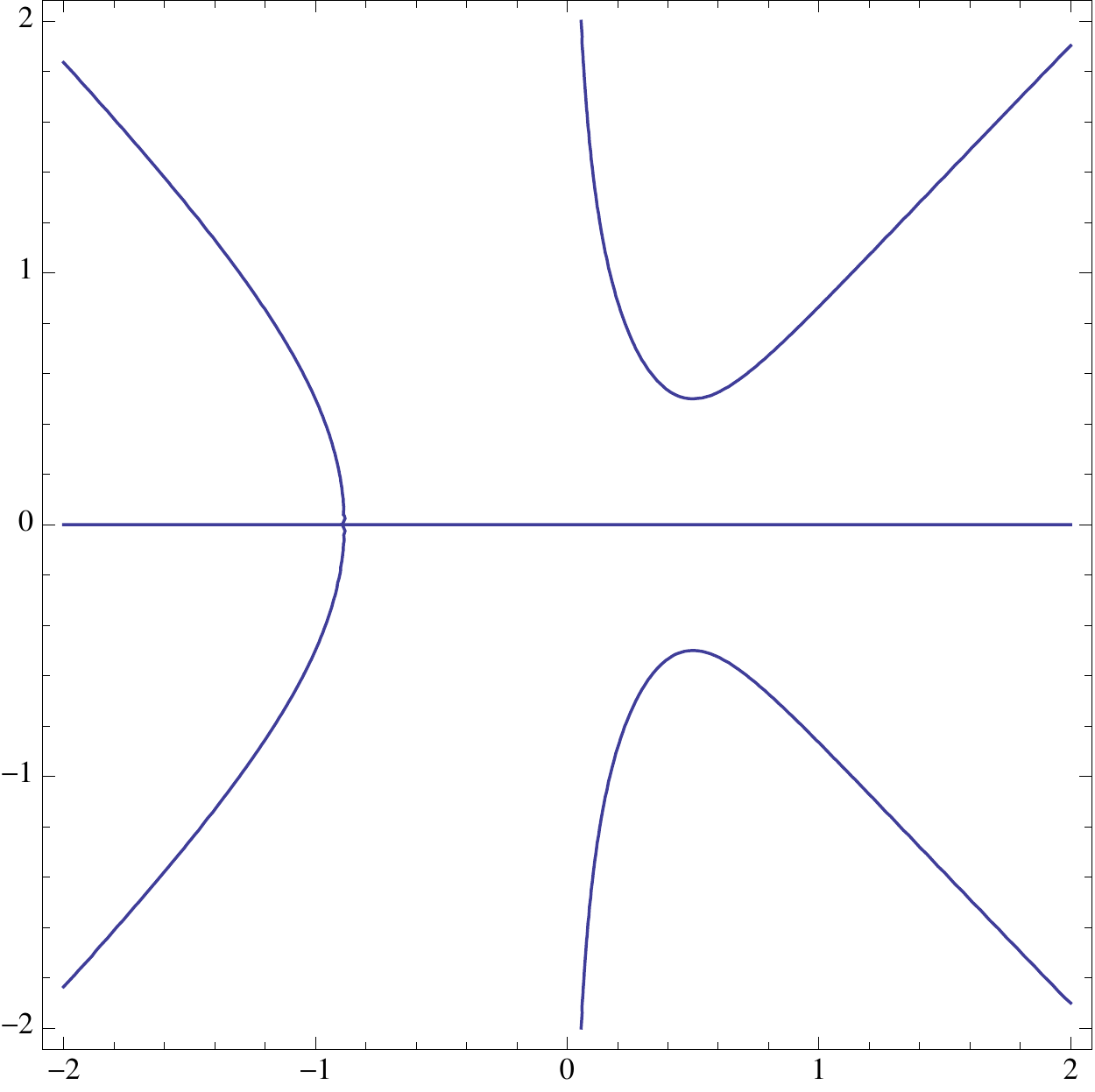}
      \put(15,-13){$(a)$ $\omega > -3(\alpha c_2)^{2/3}$}
\end{overpic}
\quad
\begin{overpic}[width=.3\textwidth]{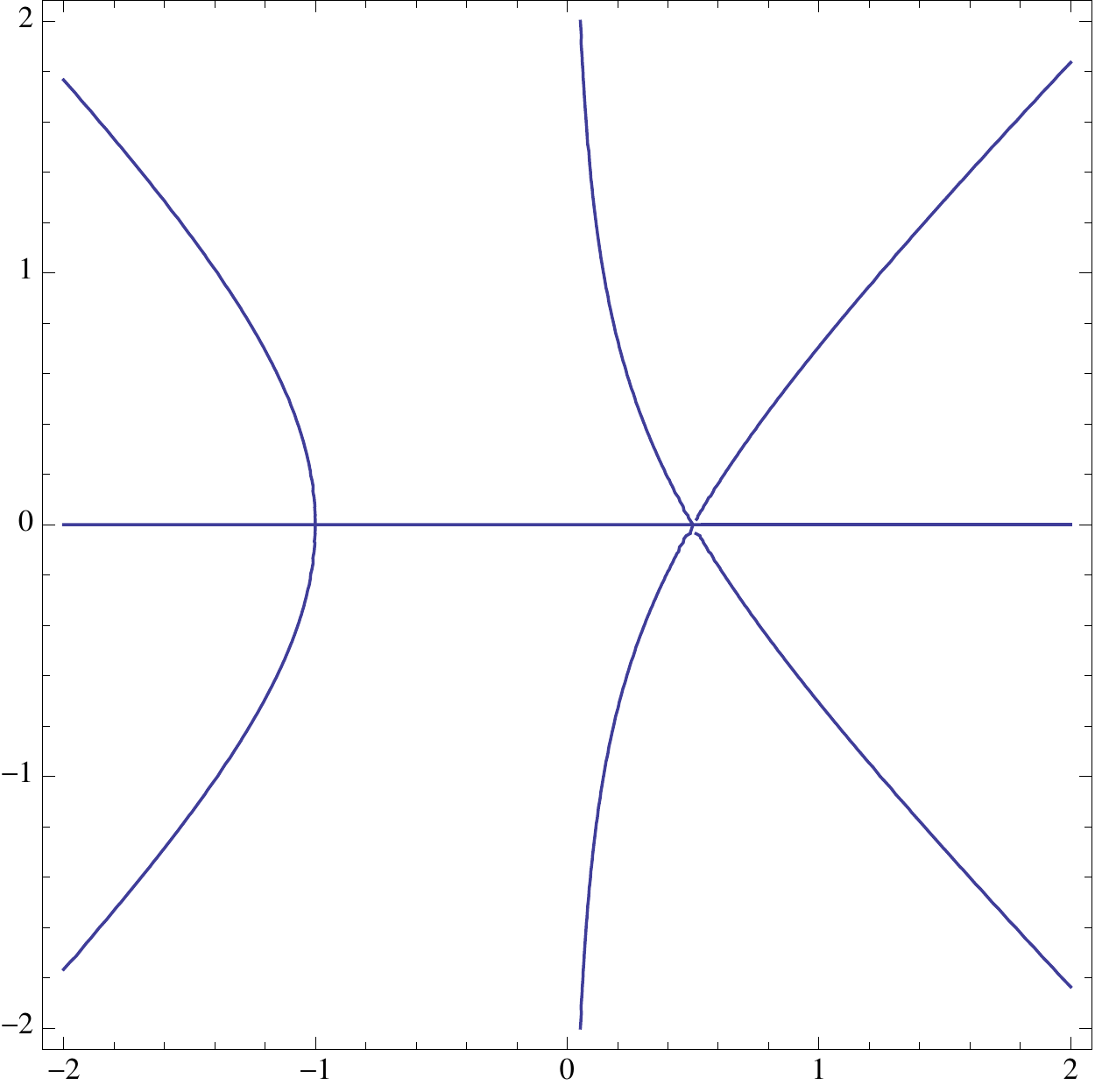}
      \put(15,-13){$(b)$ $\omega = -3(\alpha c_2)^{2/3}$}
\end{overpic}
\quad
\begin{overpic}[width=.3\textwidth]{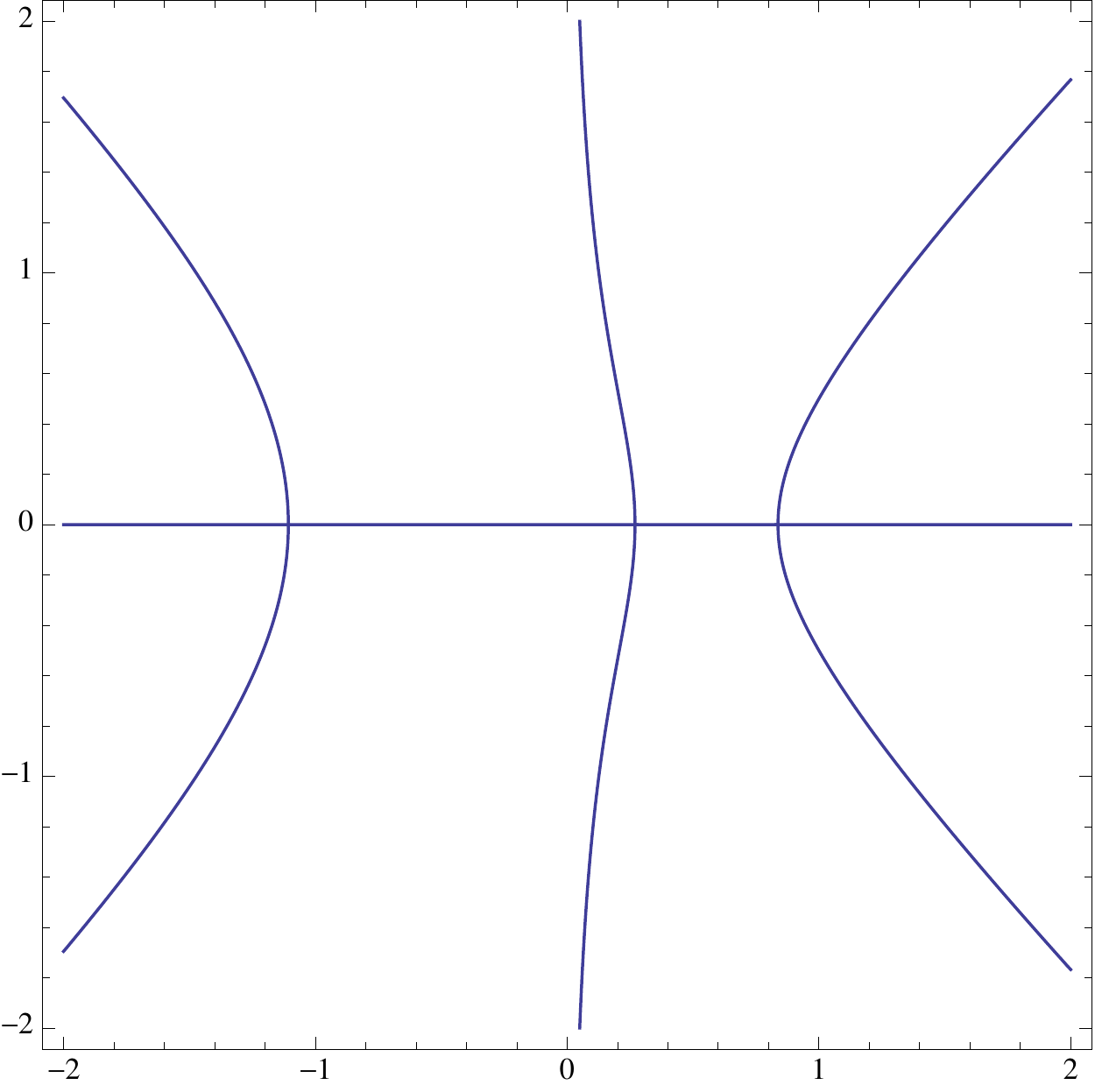}
      \put(15,-13){$(c)$ $\omega < -3(\alpha c_2)^{2/3}$}
\end{overpic}
\vspace{.7cm}
     \begin{figuretext}\label{implus}
       The set $\{\im \Omega^2(k) = 0\} = \R \cup \Gamma$ in the case of $c_2> 0$ and different values of $\omega \in \R$.
     \end{figuretext}
     \end{center}
\end{figure}

\begin{figure}
\begin{center}
\begin{overpic}[width=.45\textwidth]{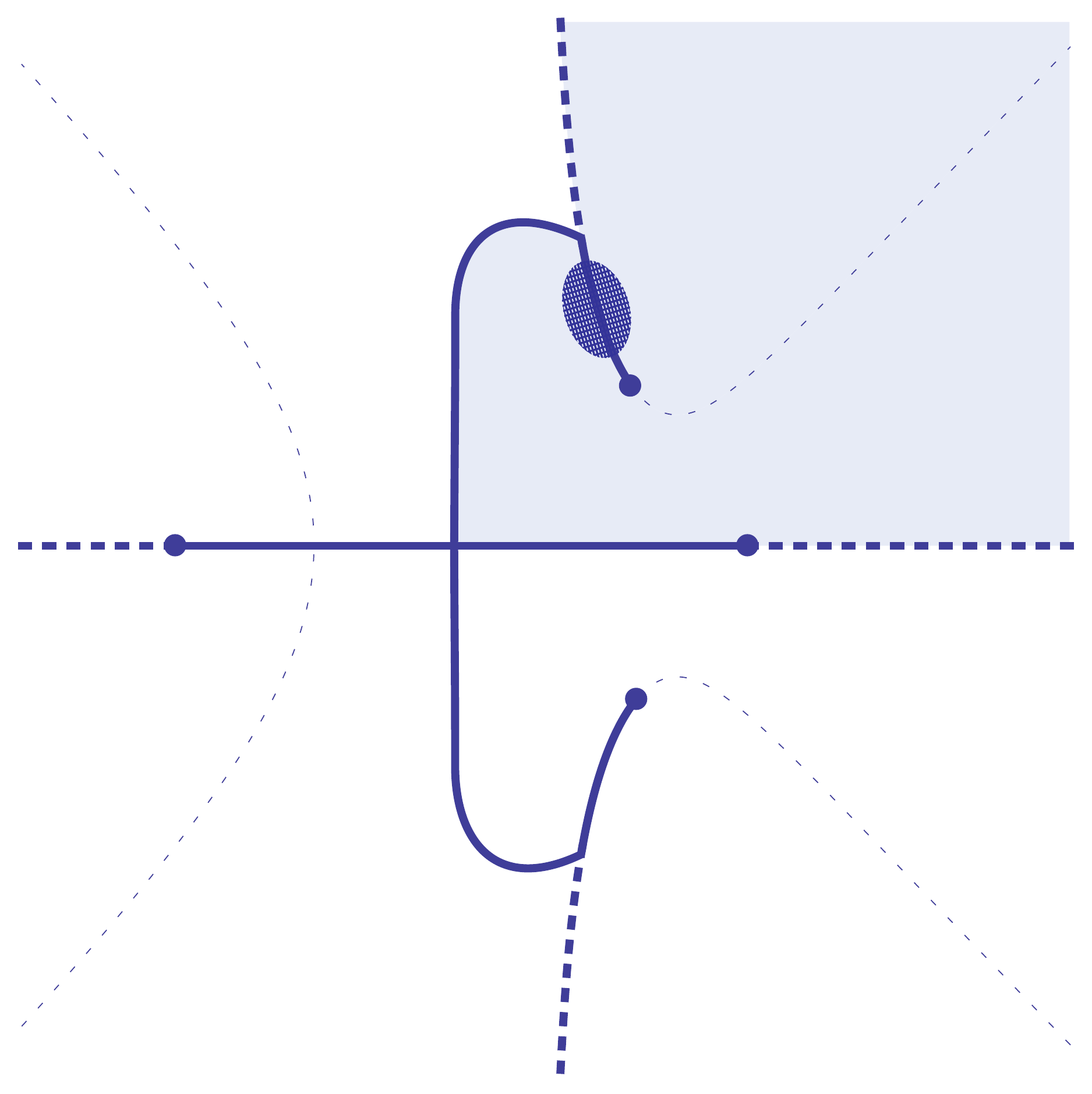}
      \put(65,75){$D_1$}
\end{overpic}
     \begin{figuretext}\label{cuts1.pdf}
        An inadmissible situation in the case of $c_2 > 0$ and $\omega > -3|\alpha c_2|^{2/3}$.
         In this and other similar figures in this section, the light shaded area is $D_1$. The dark shaded area is $U$. The solid thick lines are the branch cuts connecting the four zeros (indicated by dots) of $\Omega^2(k)$. Furthermore, $\im \Omega(k) = 0$ on the thick striped lines, while the thin striped lines indicate the set $\{\im \Omega^2(k) = 0\} \cap \{\re \Omega^2(k) < 0\}$.
         The branch cut connecting $\frac{i\bar{c}}{2\alpha}$ and $-\frac{ic}{2\alpha}$ is not displayed.
     \end{figuretext}
     \end{center}
\end{figure}

The qualitative structure of the set $\{\im \Omega^2(k) = 0\}$ is determined by the values of $c_2$ and $\omega$. We will prove Theorem \ref{mainth} by considering in turn the three cases (i) $c_2 > 0$, (ii) $c_2 = 0$, and (iii) $c_2 < 0$ including their various subcases. We will use Lemma \ref{inadmissiblelemma} to rule out all triples $(\alpha,\omega, c)$ that are not included in (\ref{admissiblesets}) as inadmissible.  

\subsection{$c_2 > 0$}
If $c_2 > 0$, the set $\{\im \Omega^2(k) = 0\}$ has the qualitative structure shown in Figure \ref{implus}. Thus, there are three subcases depending on whether $\omega$ is larger than, equal to, or smaller than $-3(\alpha c_2)^{2/3}$.

\begin{figure}
\begin{center}
\begin{overpic}[width=.45\textwidth]{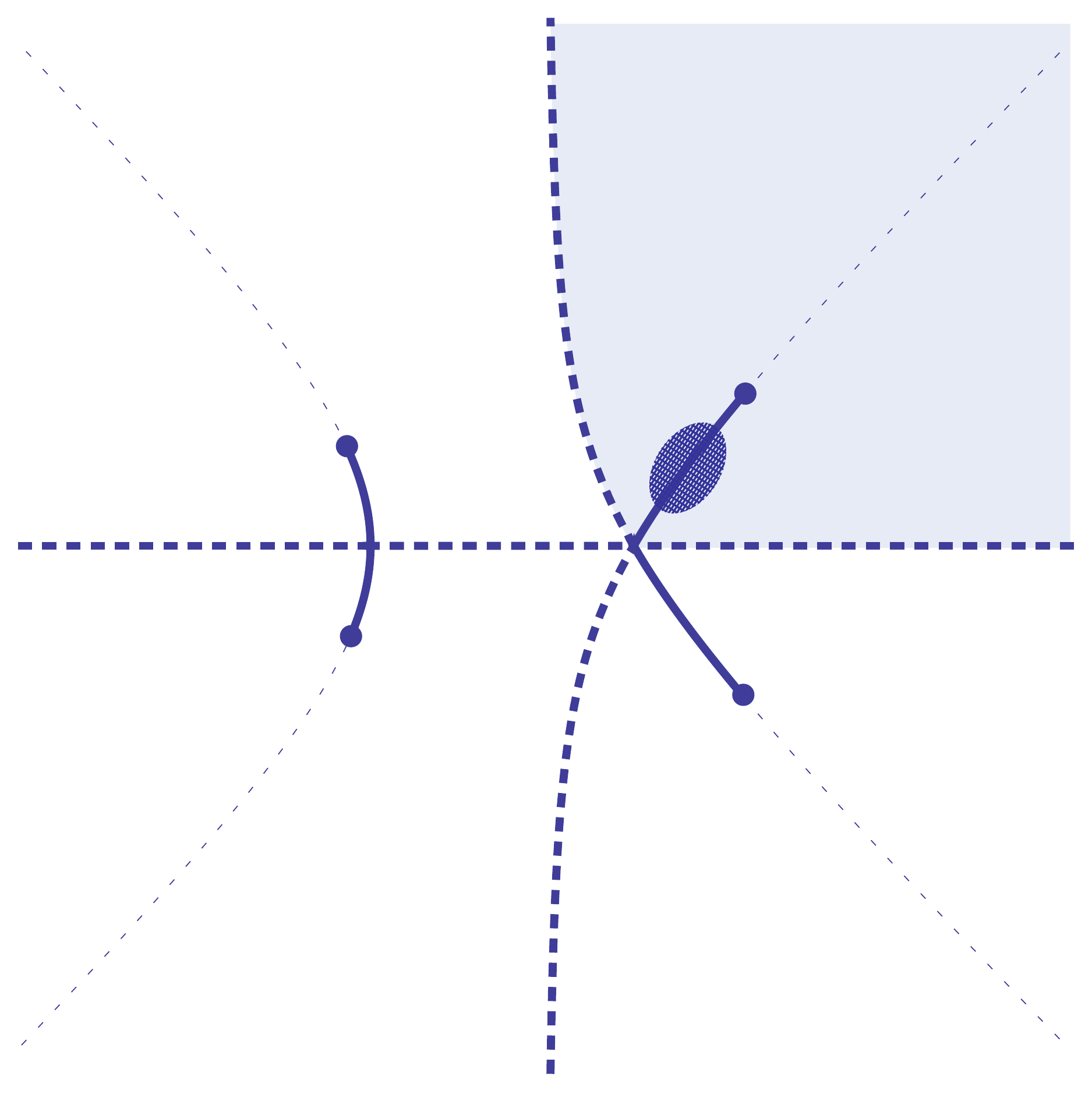}
      \put(65,75){$D_1$}
\end{overpic}
     \begin{figuretext}\label{cuts2a.pdf}
    An inadmissible situation in the case of $c_2 > 0$ and $\omega = -3(\alpha c_2)^{2/3}$.
       \end{figuretext}
     \end{center}
\end{figure}

\begin{figure}
\begin{center}
\begin{overpic}[width=.45\textwidth]{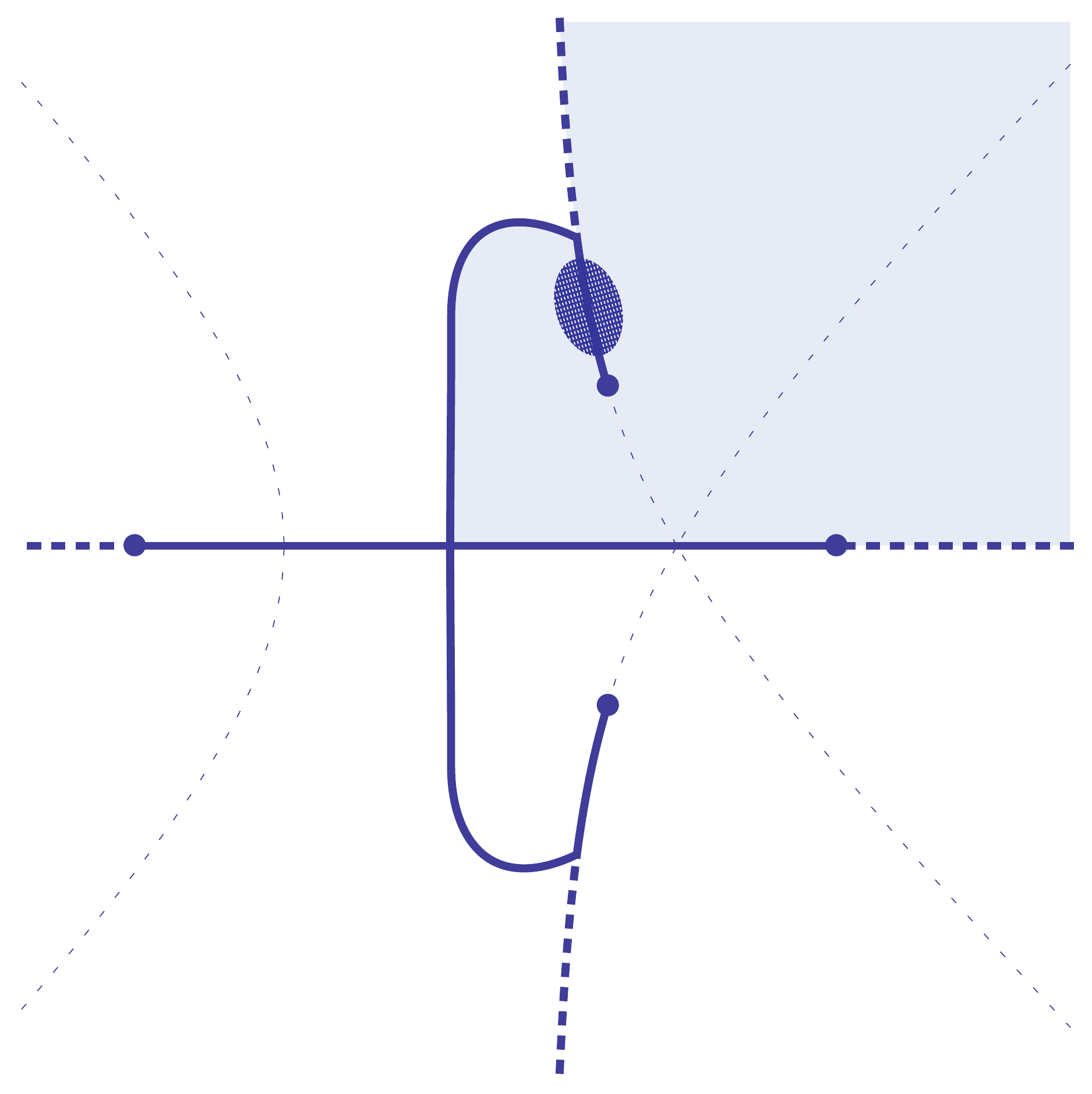}
      \put(65,75){$D_1$}
\end{overpic}
     \begin{figuretext}\label{cuts2b.pdf}
     Another inadmissible situation in the case of $c_2 > 0$ and $\omega = -3(\alpha c_2)^{2/3}$.
      \end{figuretext}
     \end{center}
\end{figure}

\begin{figure}
\begin{center}
\begin{overpic}[width=.45\textwidth]{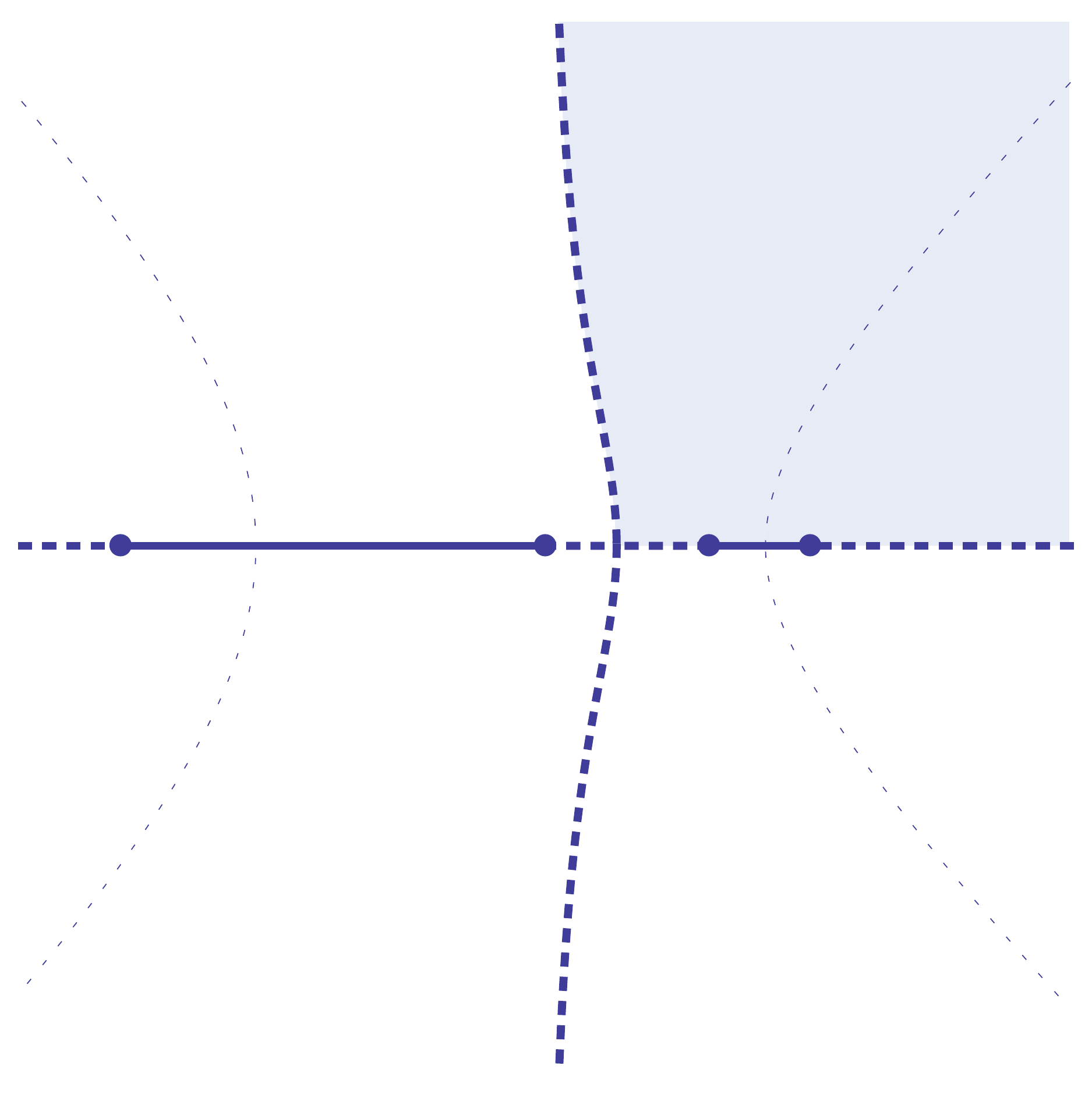}
      \put(65,75){$D_1$}
\end{overpic}
     \begin{figuretext}\label{cuts3a.pdf}
     An inadmissible situation in the case of $c_2 > 0$ and $\omega < -3(\alpha c_2)^{2/3}$.
      \end{figuretext}
     \end{center}
\end{figure}

\begin{figure}
\begin{center}
\begin{overpic}[width=.45\textwidth]{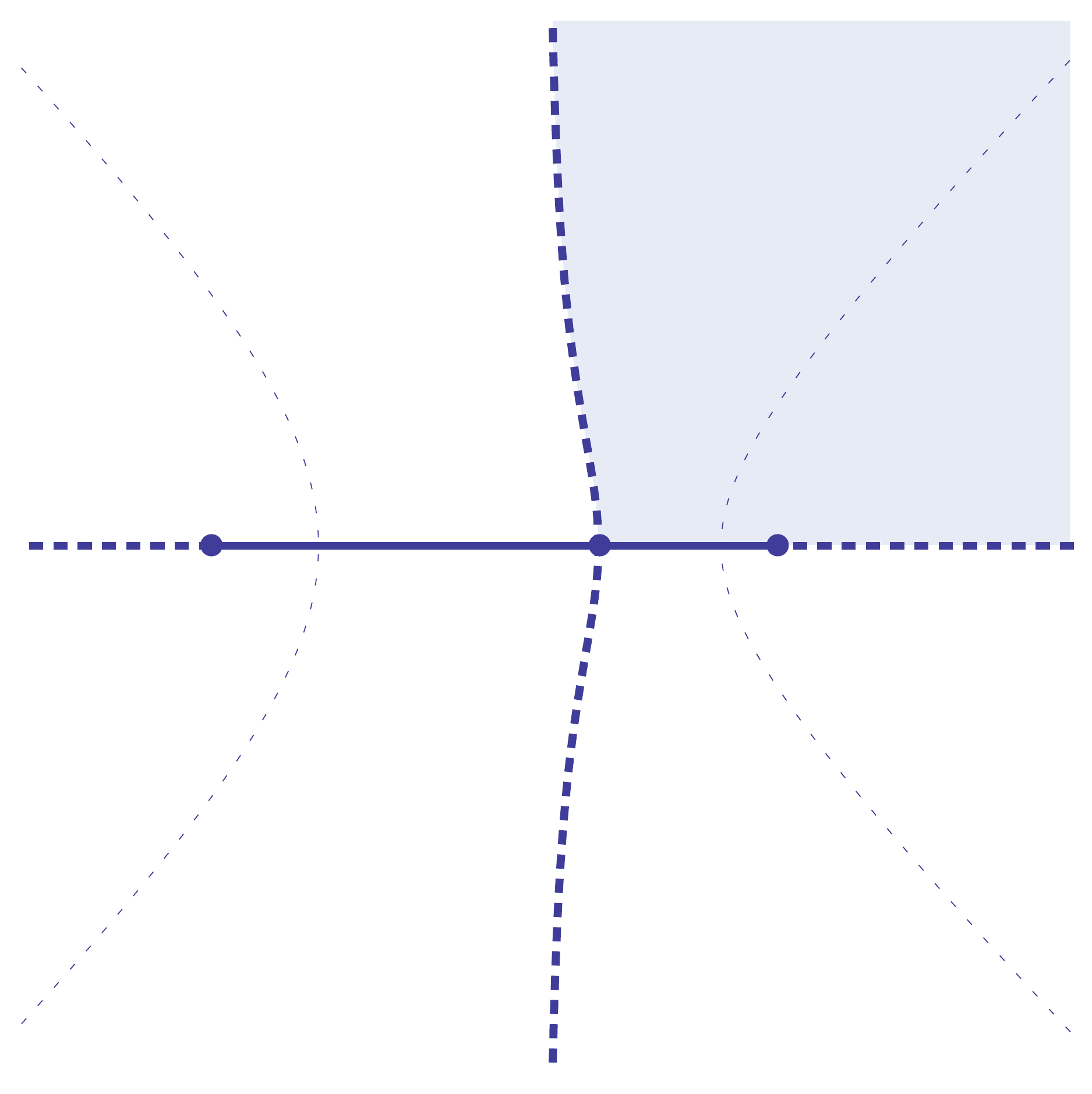}
      \put(65,75){$D_1$}
\end{overpic}
     \begin{figuretext}\label{cuts3b.pdf}
       Another inadmissible situation in the case of $c_2 > 0$ and $\omega < -3(\alpha c_2)^{2/3}$.
      \end{figuretext}
     \end{center}
\end{figure}

\subsubsection{$c_2 > 0$, $\omega > -3(\alpha c_2)^{2/3}$}
In this case, there exists a pair of zeros $\{K, \bar{K}\}$ of $\Omega^2(k)$ such that $\re K > 0$ and $\im K > 0$; the other two zeros of $\Omega^2(k)$ are real or lie in the left half-plane. Regardless of the position of the other two zeros, by choosing the branch cut connecting $K$ with $\bar{K}$ and the neighborhood $U$ as in Figure \ref{cuts1.pdf}, we infer that all these triples are inadmissible.

\subsubsection{$c_2 > 0$, $\omega = -3(\alpha c_2)^{2/3}$}
In this case, the set $\Gamma$ intersects the real axis at $\frac{(\alpha c_2)^{1/3}}{2} = \sqrt{\frac{|\omega|}{12}}$ and at $-(\alpha c_2)^{1/3} = -\sqrt{\frac{|\omega|}{3}}$. 
Thus, there are two possibilities: (i) there exists a pair of zeros $\{K, \bar{K}\}$ of $\Omega^2(k)$ such that $\re K > 0$ and $\im K > 0$, or (ii) $K = \sqrt{\frac{|\omega|}{12}}$ is a zero of $\Omega^2(k)$.
Case (i) leads to inadmissible triples, see Figures \ref{cuts2a.pdf} and \ref{cuts2b.pdf}.

Letting $\omega = -3(\alpha c_2)^{2/3}$ in the definition (\ref{Omega2def}) of $\Omega^2(k)$, we see that case (ii) occurs iff
$$c_1^2 = \frac{\left(3 \alpha^2+\omega \right)^3}{27 \alpha^2},$$ 
and if this is the case, then $\sqrt{\frac{|\omega|}{12}}$ is necessarily a triple zero. We cannot use Lemma \ref{inadmissiblelemma} to rule out these triples. Thus case (ii)  gives rise to the two-parameter family (\ref{admissibleA}) of potentially admissible triples (see Figure \ref{cutsA.pdf} for the structure of the set $\im \Omega = 0$ in this case).

\subsubsection{$c_2 > 0$, $\omega < -3(\alpha c_2)^{2/3}$}
In this case, the set $\Gamma$ intersects the real axis in three distinct points $\{K_j\}_1^3$ which we label so that $K_1 < K_2 < K_3$. 
If there exists a pair of zeros $\{K, \bar{K}\}$ of $\Omega^2(k)$ such that $\re K > 0$ and $\im K > 0$, then the triple is inadmissible as a consequence of Lemma \ref{inadmissiblelemma}. 

Let us consider the case where $\Omega^2$ has at least one simple real zero $K$ with $K > K_2$, see Figure \ref{cuts3a.pdf}. It appears difficult to rule out these triples using Lemma \ref{inadmissiblelemma}. However, these triples are inadmissible as a consequence of an argument given in \cite{FLtperiodic}, which shows that $\Omega^2$ cannot have any zeros of odd order in $\text{Int}(\bar{D}_1 \cup \bar{D}_4) \cap \R$ (see Section 3 of \cite{FLtperiodic}). 

It only remains to consider the cases where $k = K_3$ is a double zero of $\Omega^2$ and the other two zeros are real $\leq K_2$ or a complex conjugate pair in the left-half plane.  
According to (\ref{doubleconditions}), $\Omega^2(k)$ has a zero of order $\geq 2$ at $K > 0$ for some $c_1 \geq 0$ iff
\begin{align}\label{4K3omegaK}
4K^3 + \omega K + \alpha c_2 = 0 \quad \text{and} \quad \bigg(\alpha^2+ \frac{\omega}{2}\bigg)^2 -c_2^2 - 2K^2(6K^2 + \omega) \geq 0.
\end{align}
In the present case where $\alpha >0$, $c_2 > 0$, and $\omega < -3(\alpha c_2)^{2/3}$, this happens iff either
\begin{align}\label{case31}
\omega < -12 K^2, \quad 0 < c_2 \leq -\frac{4K^2 + \omega}{2}, \quad \alpha = -\frac{4K^3 + \omega K}{c_2},
\end{align}
or 
\begin{align}\label{case32}
-12 K^2 < \omega < -4K^2, \quad 0 < c_2 \leq -\frac{4K^2 + \omega}{2}, \quad \alpha = -\frac{4K^3 + \omega K}{c_2},
\end{align}
or 
\begin{align}\label{case33}
-12 K^2 < \omega < -4K^2, \quad c_2 = \sqrt{-8K^4 - 2K^2 \omega}, \quad \alpha = -\frac{4K^3 + \omega K}{c_2}.
\end{align}

In the case of (\ref{case31}), since the first equation in (\ref{4K3omegaK}) implies that $K = K_j$ for some $j$ and $\frac{d}{dK}(4K^3 + \omega K + \alpha c_2) < 0$, the double zero $K$ coincides with the intersection point $K_2$, see Figure \ref{cuts3b.pdf}. Since $\Omega^2$ then has a simple real zero strictly bigger than $K_2$, this case is inadmissible by the argument of \cite{FLtperiodic}.

In the cases (\ref{case32}) and (\ref{case33}) the double zero $K$ coincides with $K_3$. Hence these cases cannot be ruled out by means of Lemma \ref{inadmissiblelemma}. 
The case (\ref{case32}) gives rise to the family (\ref{admissibleB}) (see Figure \ref{cutsB.pdf}).
In the case of (\ref{case33}), we have
$$c_1^2 = -c_2^2 - 2K^2(6K^2 + \omega) + \bigg(\alpha^2+ \frac{\omega}{2}\bigg)^2 = 0.$$
Hence, using $\alpha > 0$ as a parameter instead of $K > 0$, we find that (\ref{case33}) gives rise to the family (\ref{admissibleC}) of potentially admissible triples (see Figure \ref{cutsC.pdf}).

%%%%%%%%%%% im c = 0 %%%%%%%%%%%%%%%%%%%%%%%%%%%%%%%
\subsection{$c_2 = 0$} 
\begin{figure}
\begin{center}
\begin{overpic}[width=.3\textwidth]{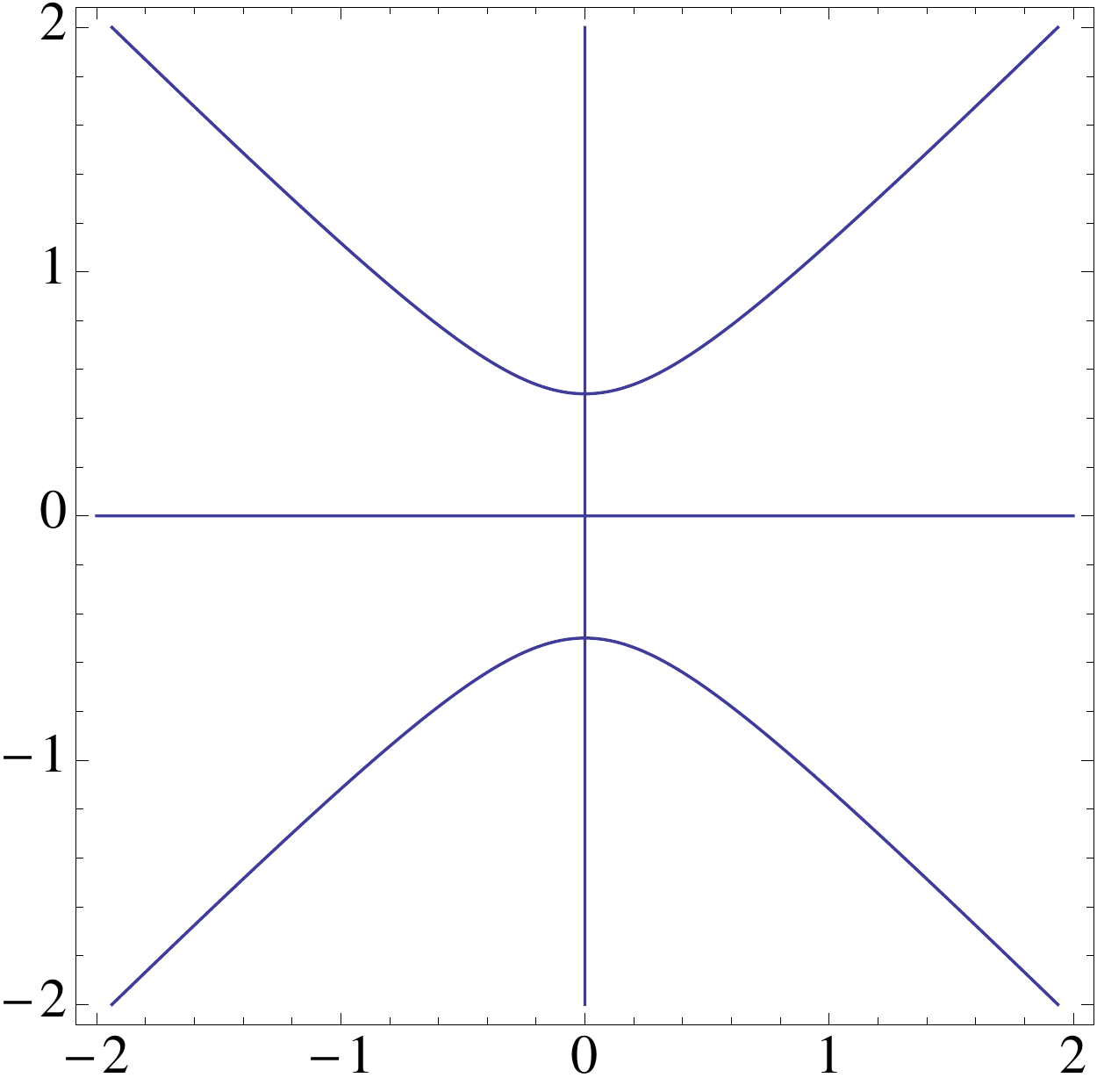}
      \put(15,-13){$(a)$ $\omega > 0$}
\end{overpic}
\quad
\begin{overpic}[width=.3\textwidth]{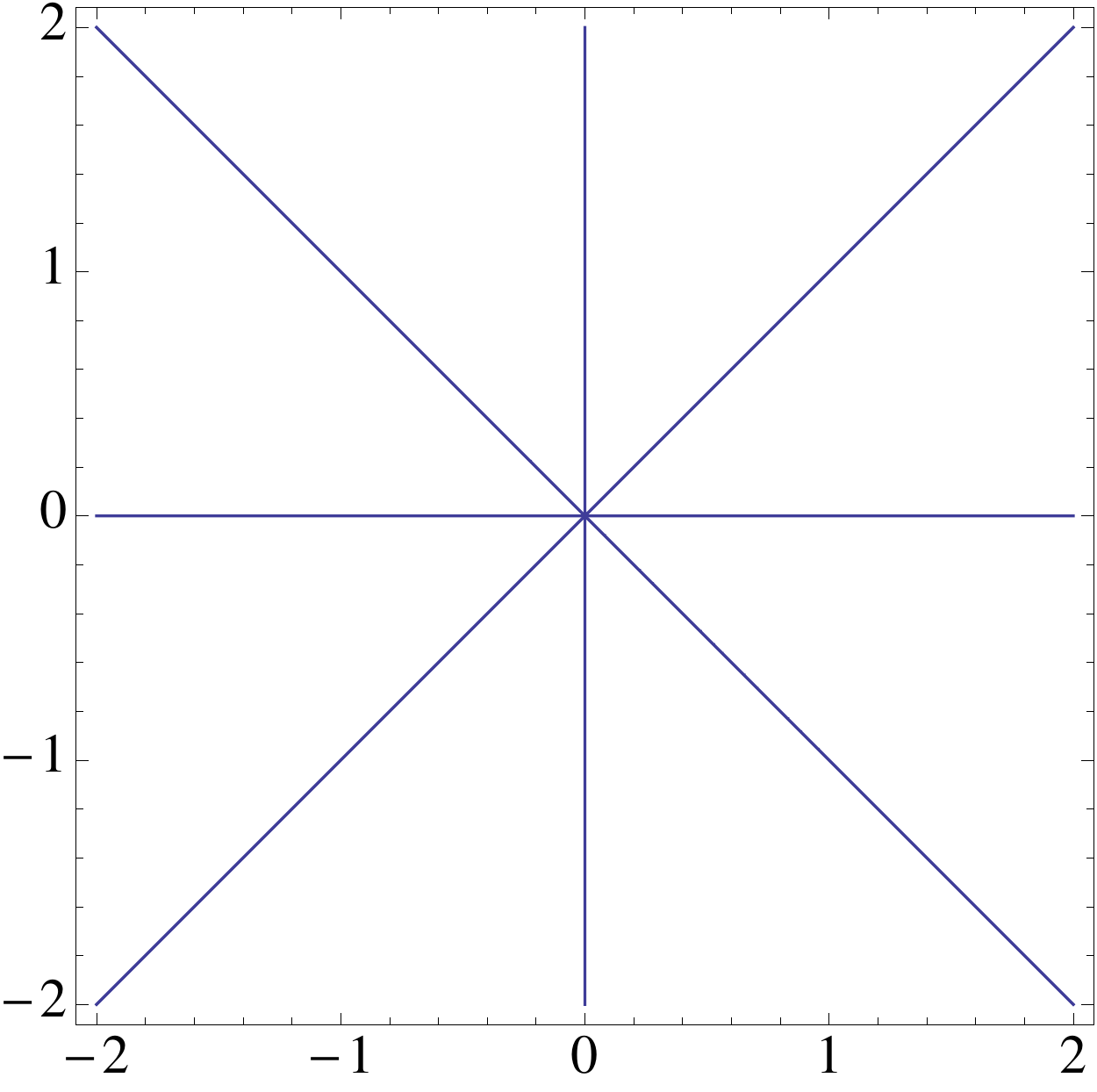}
      \put(15,-13){$(b)$ $\omega = 0$}
\end{overpic}
\quad
\begin{overpic}[width=.3\textwidth]{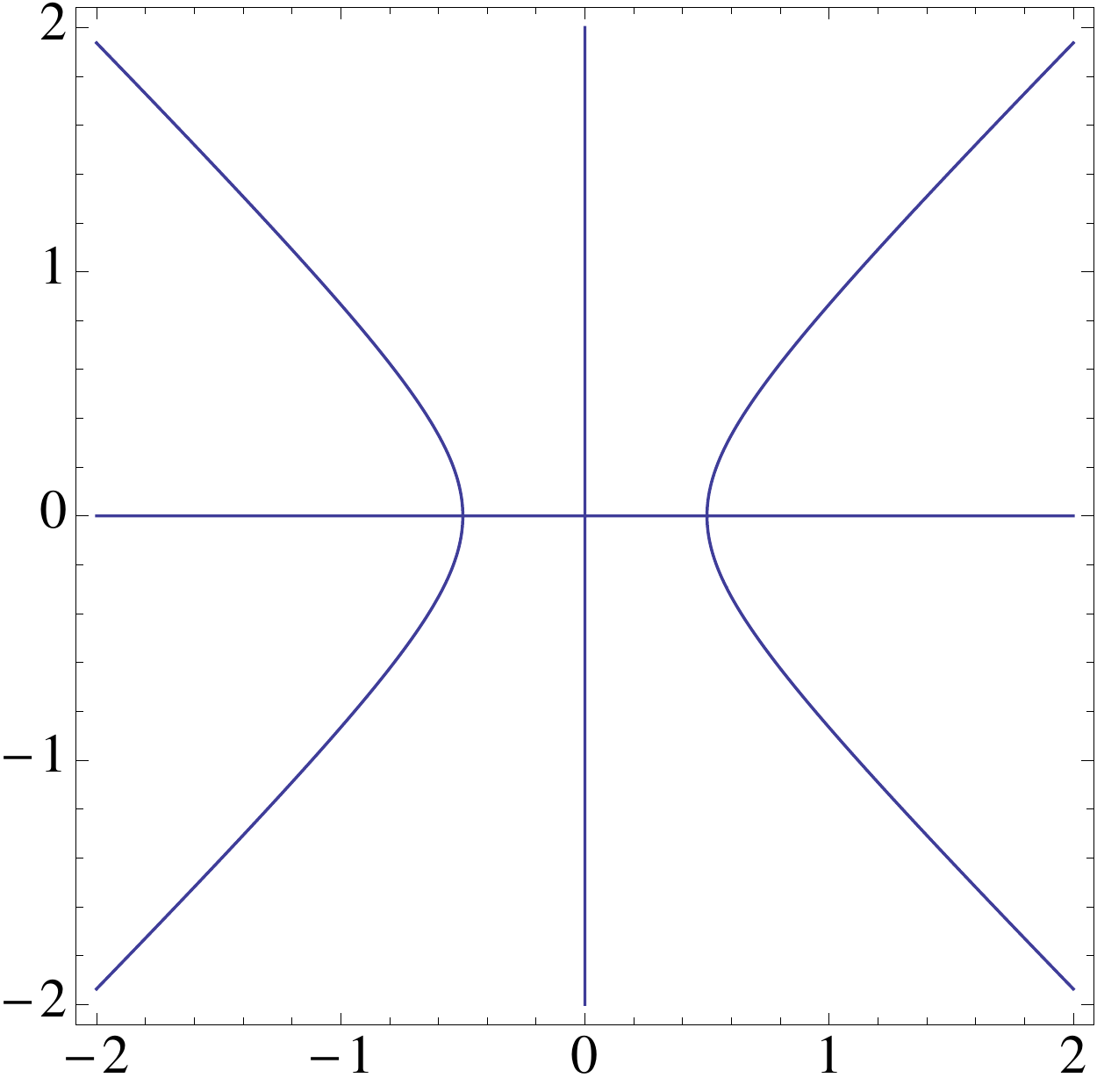}
      \put(15,-13){$(c)$ $\omega < 0$}
\end{overpic}
\vspace{.7cm}
     \begin{figuretext}\label{imzero}
       The set $\{\im \Omega^2(k) = 0\}  = \R \cup \Gamma$ in the case of $c_2=0$ and different values of $\omega \in \R$.
     \end{figuretext}
     \end{center}
\end{figure}

The structure of the set $\{\im \Omega^2(k) = 0\}$ is shown in Figure \ref{imzero}.
The zero set of $\Omega^2(k)$ is invariant under each of the two reflections $k_2 \to -k_2$ and $k_1 \to -k_1$ in the real and imaginary axes, respectively.
 
\subsubsection{$c_2 = 0$, $\omega > 0$}
The set $\{\im \Omega^2(k) = 0\}$ consists of the real and imaginary axes together with the two parabolas
$$k_2 = \pm \sqrt{k_1^2 + \frac{\omega}{4}},$$
which intersect the imaginary axis at $k = \pm \frac{i\sqrt{\omega}}{2}$. 

\begin{figure}
\begin{center}
\begin{overpic}[width=.45\textwidth]{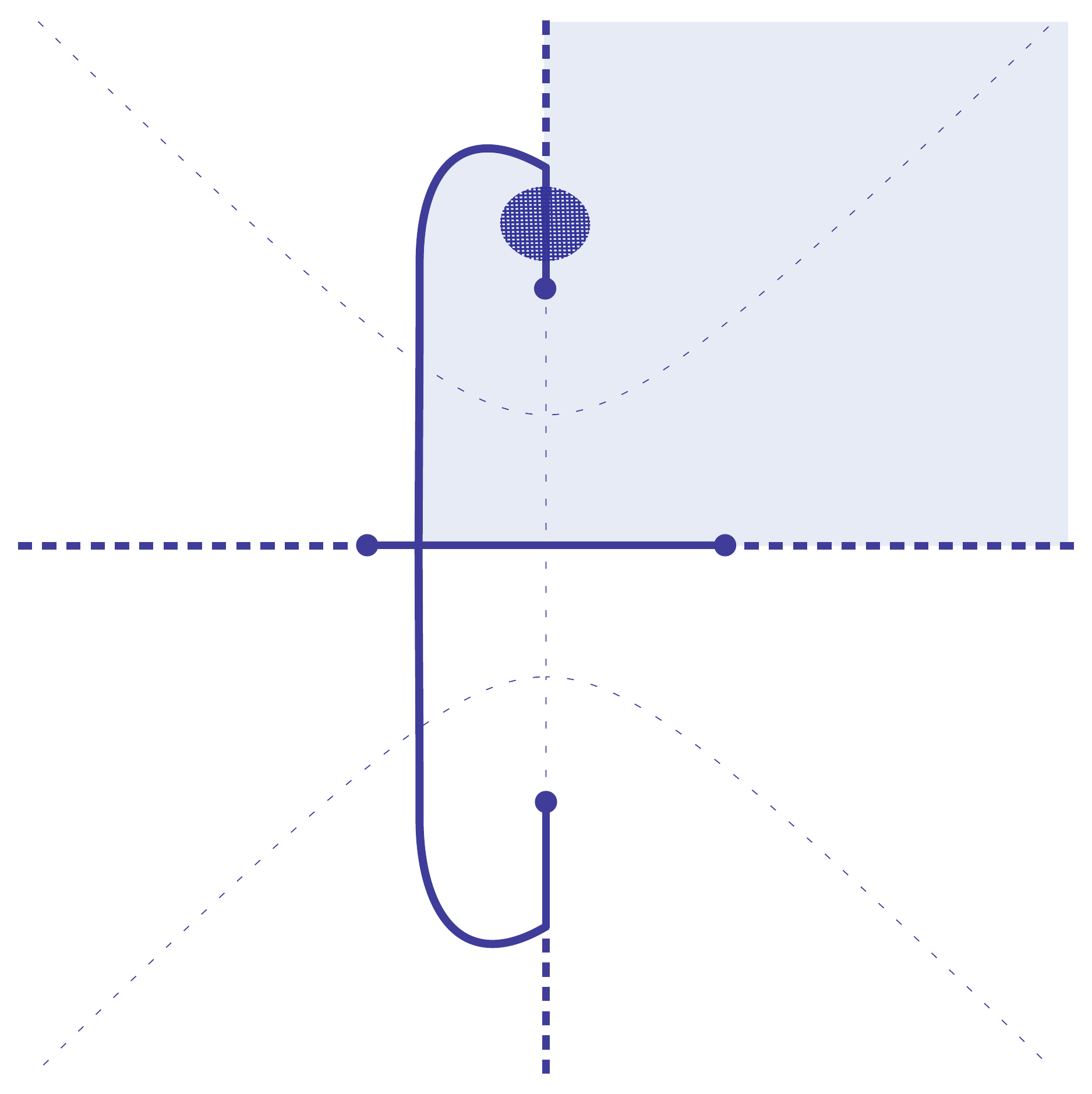}
      \put(65,82){$D_1$}
\end{overpic}
     \begin{figuretext}\label{cuts4a.pdf}
     An inadmissible situation in the case of $c_2 = 0$ and $\omega > 0$.
      \end{figuretext}
     \end{center}
\end{figure}

\begin{figure}
\begin{center}
\begin{overpic}[width=.45\textwidth]{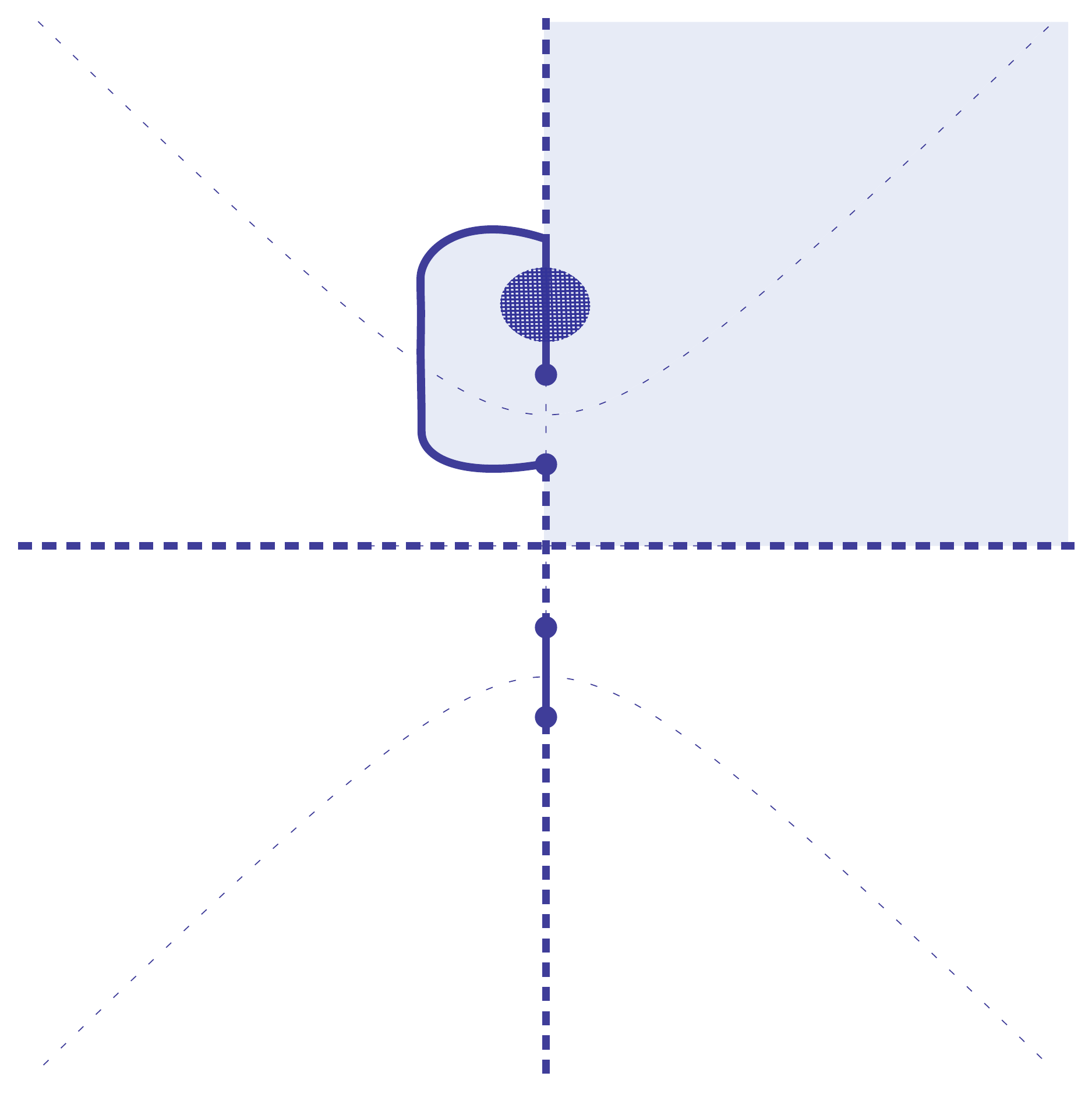}
      \put(65,82){$D_1$}
\end{overpic}
     \begin{figuretext}\label{cuts4b.pdf}
     Another inadmissible situation in the case of $c_2 = 0$ and $\omega > 0$.
      \end{figuretext}
     \end{center}
\end{figure}

Proceeding as in the case of $c_2 > 0$, we find that $(\alpha, \omega, c)$ is inadmissible unless $\pm \frac{i\sqrt{\omega}}{2}$ are double zeros of $\Omega^2(k)$; two examples of inadmissible situations are shown in Figures \ref{cuts4a.pdf} and \ref{cuts4b.pdf}. 

Now $\pm \frac{i\sqrt{\omega}}{2}$ are double zeros of $\Omega^2(k)$ iff
$$4k^4 + 2\omega k^2 + \bigg(\frac{\omega}{2} + \alpha^2\bigg)^2 - c_1^2
= 4\bigg(k - \frac{i\sqrt{\omega}}{2}\bigg)^2\bigg(k + \frac{i\sqrt{\omega}}{2}\bigg)^2,$$
that is, iff
$$c_1^2 = \alpha^2(\omega + \alpha^2).$$
This yields the following family of potentially admissible triples which is included in (\ref{admissibleD}) (see Figure \ref{cutsDa.pdf}):
\begin{align}\label{familyDatriples}
\{(\alpha,\omega,c = \pm \alpha\sqrt{\omega + \alpha^2}) \; |\; \omega > 0, \; \alpha > 0\}.
\end{align}

\begin{figure}
\begin{center}
\begin{overpic}[width=.45\textwidth]{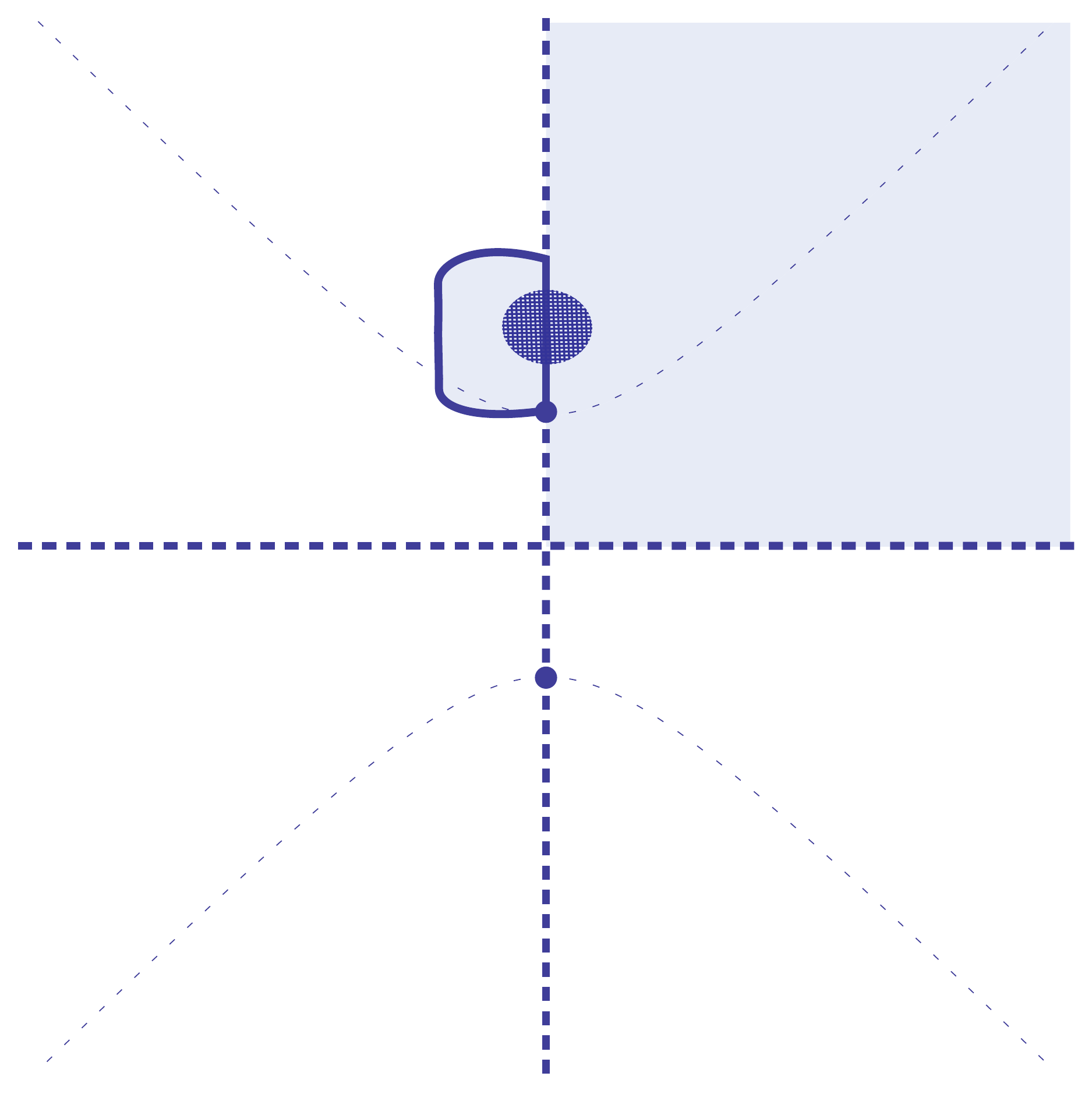}
      \put(65,82){$D_1$}
\end{overpic}
     \begin{figuretext}\label{cuts4c.pdf}
      An example of a choice of branch cuts for which $D_1 \setminus \mathcal{C}$ is not connected. Lemma \ref{inadmissiblelemma} does not apply to this situation.
      \end{figuretext}
     \end{center}
\end{figure}

\begin{remark}\upshape
One may ask why the triples (\ref{familyDatriples}) are not rendered inadmissible by the choice of branch cuts displayed in Figure \ref{cuts4c.pdf}. Figure \ref{cuts4c.pdf} contains a branch cut which starts and ends at the double zero $\frac{i\sqrt{\omega}}{2}$. 
The reason Lemma \ref{inadmissiblelemma} does not apply to the situation in Figure \ref{cuts4c.pdf} is that $D_1 \setminus \mathcal{C}$ is not connected. In fact, since $\Omega^2(k) = \frac{1}{4} (4 k^2+\omega)^2$, we have $\Omega(k) = - 2k^2- \frac{\omega}{2}$ for all $k$ surrounded by the branch cut in Figure \ref{cuts4c.pdf}. Consequently, $\im(\Omega(k) + 2k^2) = 0$ and the argument leading to the global relation (\ref{GR}) breaks down for these $k$. 
\end{remark}

\subsubsection{$c_2 = 0$, $\omega = 0$}
The triple $(\alpha, \omega, c)$ is inadmissible unless $k = 0$ is a fourth-order zero $\Omega^2(k)$, which happens iff $c_1^2 = \alpha^4$.
This yields the following family of admissible triples which is included in (\ref{admissibleD}) (see Figure \ref{cutsDb.pdf}):
$$\{(\alpha,\omega,c  = \pm \alpha\sqrt{\omega + \alpha^2}) \; | \; \omega = 0, \; \alpha > 0\}.$$

\subsubsection{$c_2 = 0$, $\omega < 0$}
The triple $(\alpha, \omega, c)$ is inadmissible unless $\pm \frac{\sqrt{|\omega|}}{2}$ are double zeros of $\Omega^2(k)$, which happens iff $c_1^2 = \alpha^2(\omega + \alpha^2)$.
This yields the following family of admissible triples which is included in (\ref{admissibleD}) (see Figure \ref{cutsDc.pdf}):
$$ \{(\alpha,\omega,c = \pm \alpha\sqrt{\omega + \alpha^2}) \; | \; \omega + \alpha^2 \geq 0,  \; \omega < 0, \; \alpha > 0\}.$$

%%%%%%%%%%% im c < 0 %%%%%%%%%%%%%%%%%%%%%%%%%%%%%%%
\subsection{$c_2 < 0$}
\begin{figure}
\begin{center}
\begin{overpic}[width=.3\textwidth]{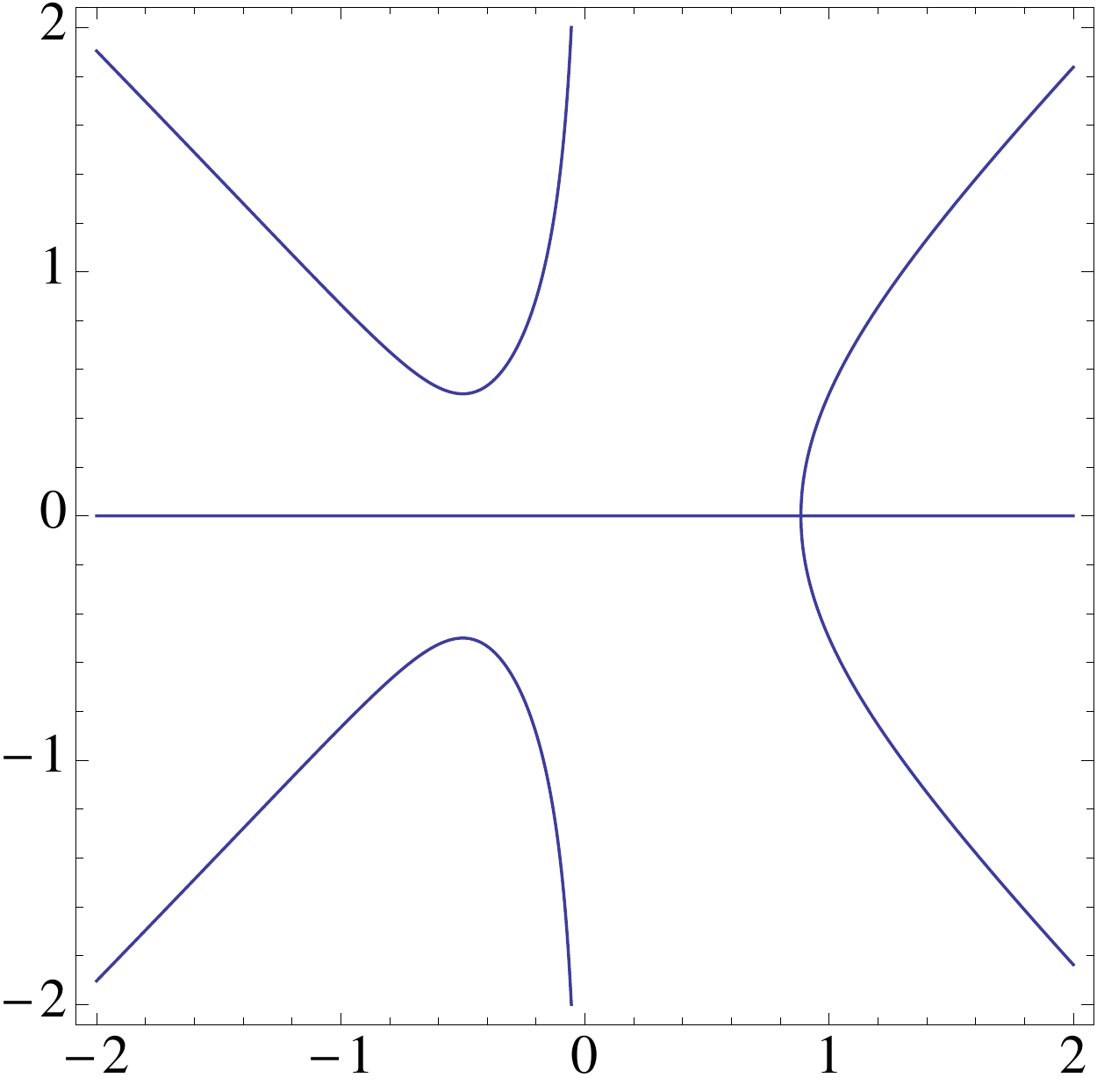}
      \put(15,-13){$(a)$ $\omega > -3|\alpha c_2|^{2/3}$}
\end{overpic}
\quad
\begin{overpic}[width=.3\textwidth]{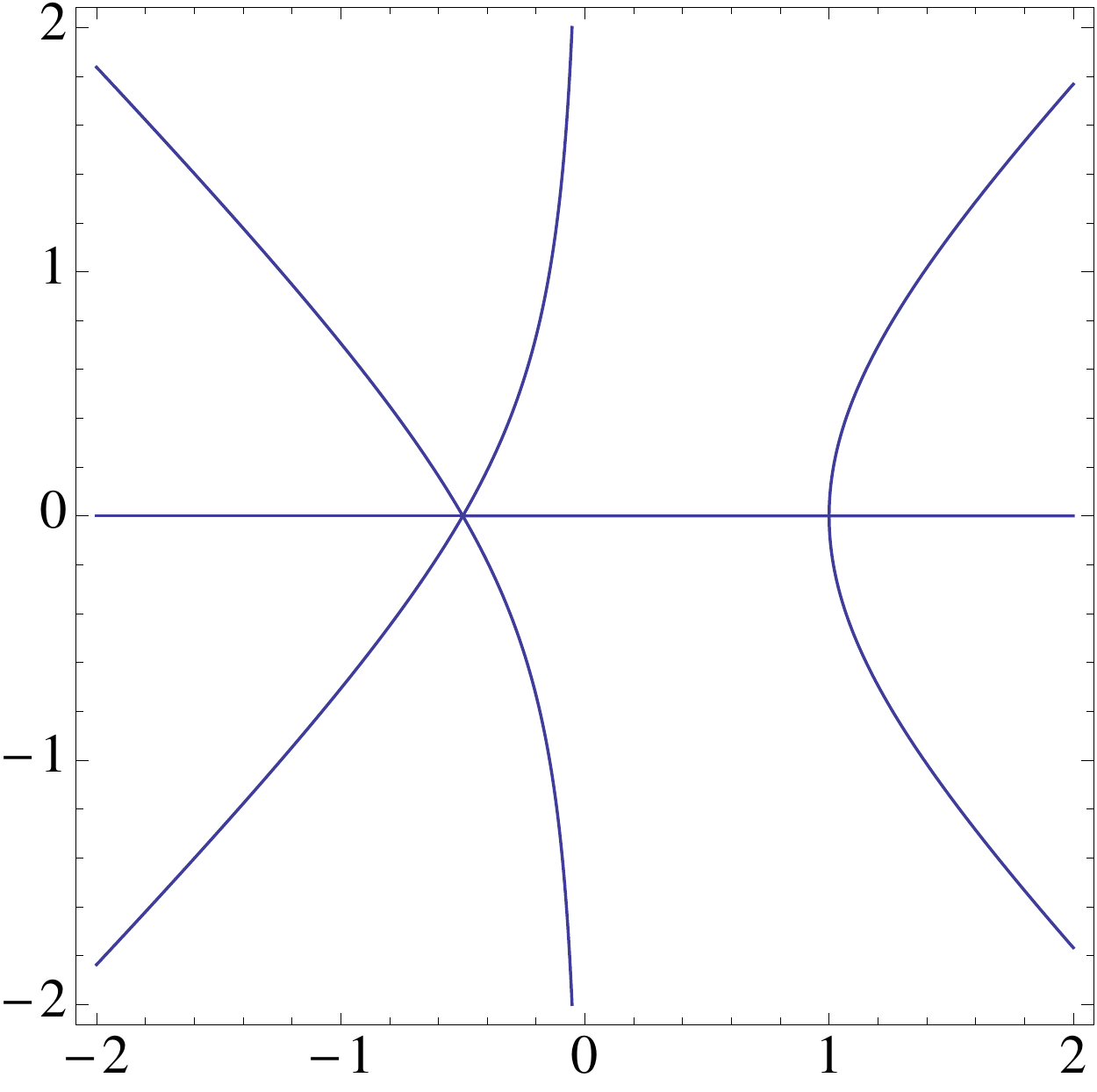}
      \put(15,-13){$(b)$ $\omega = -3|\alpha c_2|^{2/3}$}
\end{overpic}
\quad
\begin{overpic}[width=.3\textwidth]{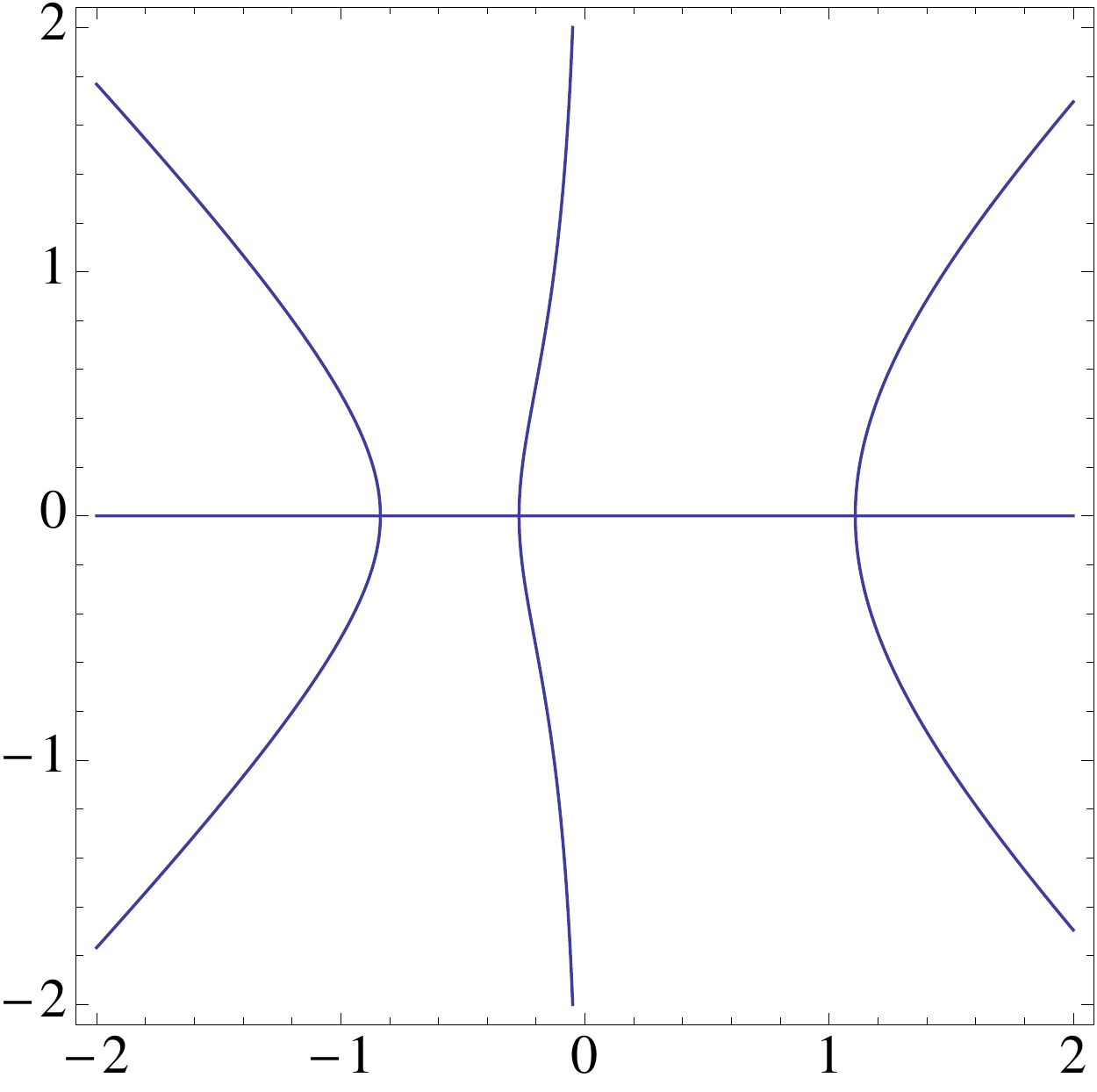}
      \put(15,-13){$(c)$ $\omega < -3|\alpha c_2|^{2/3}$}
\end{overpic}
\vspace{.7cm}
     \begin{figuretext}\label{imminus}
       The set $\{\im \Omega^2(k) = 0\}  = \R \cup \Gamma$ in the case of $c_2< 0$ and different values of $\omega \in \R$.
     \end{figuretext}
     \end{center}
\end{figure}

The structure of the set $\{\im \Omega^2(k) = 0\}$ is shown in Figure \ref{imminus}.
The sets $\{\im \Omega^2(k) = 0\}$ and $\{\re \Omega^2(k) = 0\}$ are invariant under the transformation $(c_2, k_1) \to (-c_2, -k_1)$. Hence, these sets can be obtained for $c_2 < 0$ by reflecting the corresponding sets for $c_2 > 0$ in the imaginary axis.

\subsubsection{$c_2 < 0$, $\omega > -3|\alpha c_2|^{2/3}$}
\begin{figure}
\begin{center}
\begin{overpic}[width=.45\textwidth]{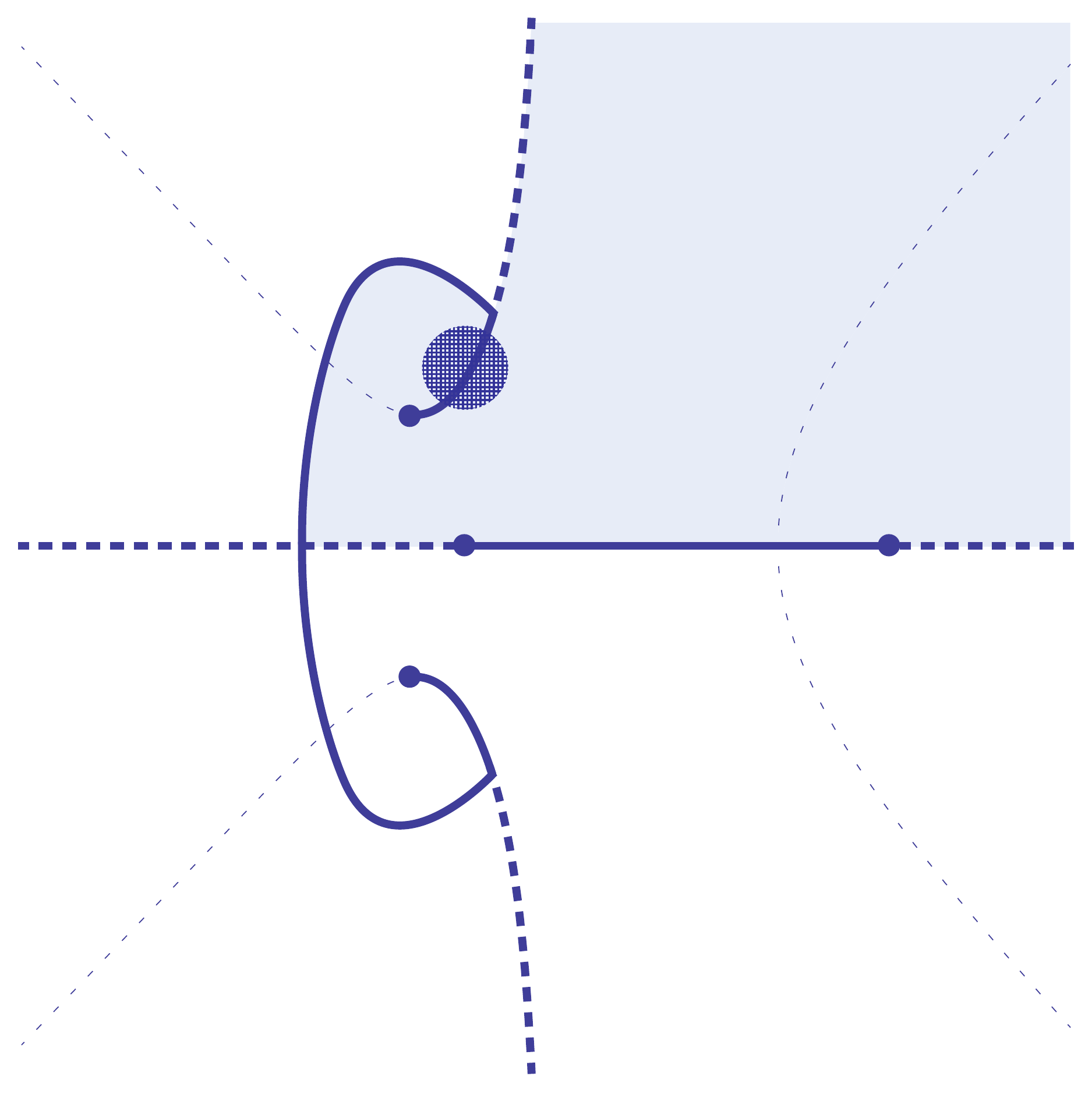}
      \put(65,75){$D_1$}
\end{overpic}
     \begin{figuretext}\label{cuts7a.pdf}
       An inadmissible situation in the case of $c_2 < 0$ and $\omega > -3|\alpha c_2|^{2/3}$.
        \end{figuretext}
     \end{center}
\end{figure}

There exists a pair of zeros $\{K, \bar{K}\}$ of $\Omega^2(k)$ such that $\re K < 0$ and $\im K > 0$, whereas the other two zeros are real or lie in the right half-plane. By choosing the branch cut connecting $K$ with $\bar{K}$ and the neighborhood $U$ as in Figure \ref{cuts7a.pdf}, we infer that all these triples are inadmissible.

\subsubsection{$c_2 < 0$, $\omega = -3|\alpha c_2|^{2/3}$}
The set $\Gamma$ intersects the real axis at $K_1 = -\frac{|\alpha c_2|^{1/3}}{2} = -\sqrt{\frac{|\omega|}{12}}$ and at $K_2 = |\alpha c_2|^{1/3} = \sqrt{\frac{|\omega|}{3}}$. 
Lemma \ref{inadmissiblelemma} implies that only triples for which $K_2$ is a double zero of $\Omega^2(k)$ can be admissible. On the other hand, if $K_2$ is a double zero then the other two zeros of $\Omega^2(k)$ have real parts $\leq K_1$.
 
According to (\ref{doubleconditions}), $\Omega^2(k)$ has a double zero at $K  = K_2 > 0$ for some $c_1 \geq 0$ if and only if the equations in (\ref{4K3omegaK}) are satisfied.
In the present case where $\alpha >0$, $c_2 < 0$, and $\omega = -3|\alpha c_2|^{2/3}$, this happens if and only if 
\begin{align*}
  -\frac{K^2}{2} \leq c_2 < 0, \quad \omega = -3K^2, \quad \alpha = -\frac{4K^3 + \omega K}{c_2}.
\end{align*}
This yields the following two-parameter family of potentially admissible triples  parametrized by $K$ and $c_2$:
\begin{align*}
\bigg\{\bigg(\alpha = -\frac{4K^3 + \omega K}{c_2}, \omega = -3K^2, c = \pm \sqrt{\bigg(\alpha^2+ \frac{\omega}{2}\bigg)^2 -c_2^2 - 2K^2(6K^2 + \omega)} + ic_2\bigg) &
	\\
 \bigg| \; -\frac{K^2}{2} \leq c_2 < 0, K > 0\bigg\}. &
\end{align*}
This family is the subset of (\ref{admissibleE}) for which $\omega = -3K^2$ (see Figure \ref{cutsEa.pdf}).

\subsubsection{$c_2 < 0$, $\omega < -3|\alpha c_2|^{2/3}$}
The set $\Gamma$ intersects the real axis at the three points $K_1 < K_2 < 0 <  K_3$. The only possibly admissible case occurs when $K_3$ is a double zero of $\Omega^2(k)$. 
According to (\ref{doubleconditions}), $\Omega^2(k)$ has a double zero at $K = K_3 > 0$ for some $c_1 \geq 0$ if and only if the equations in (\ref{4K3omegaK}) are satisfied.
In the present case where $\alpha >0$, $c_2 < 0$, and $\omega < -3|\alpha c_2|^{2/3}$, this happens iff 
\begin{align*}
  -\frac{4K^2 + \omega}{2} \leq c_2 < 0, \quad -4K^2 < \omega < -3K^2, \quad \alpha = -\frac{4K^3 + \omega K}{c_2}.
\end{align*}
This yields the following three-parameter family of potentially admissible triples parametrized by $K$, $\omega$, and $c_2$:
\begin{align*}
\bigg\{\bigg(\alpha = -\frac{4K^3 + \omega K}{c_2}, \omega, c = \pm \sqrt{\bigg(\alpha^2+ \frac{\omega}{2}\bigg)^2 -c_2^2 - 2K^2(6K^2 + \omega)} + ic_2\bigg) \qquad\qquad & 
	\\
\bigg| \; 
-4K^2 < \omega < -3K^2, \; -\frac{4K^2 + \omega}{2} \leq c_2 < 0, \; K > 0\bigg\}.& 
\end{align*}
This family is the subset of (\ref{admissibleE}) for which $-4K^2 < \omega < -3K^2$ (see Figure \ref{cutsEb.pdf}).

This completes the proof of Theorem \ref{mainth}.

\section{Analysis of potentially admissible triples}\label{analysissec}\nequation
In this section we lay the foundation for determining which of the potentially admissible triples of Theorem \ref{mainth} are actually admissible by determining the structure of the domains $\{D_j\}_1^4$ and the associated branch cuts for each of the families in (\ref{admissiblesets}).

\subsection{The family (\ref{admissibleA})}
\begin{figure}
\begin{center}
\begin{overpic}[width=.45\textwidth]{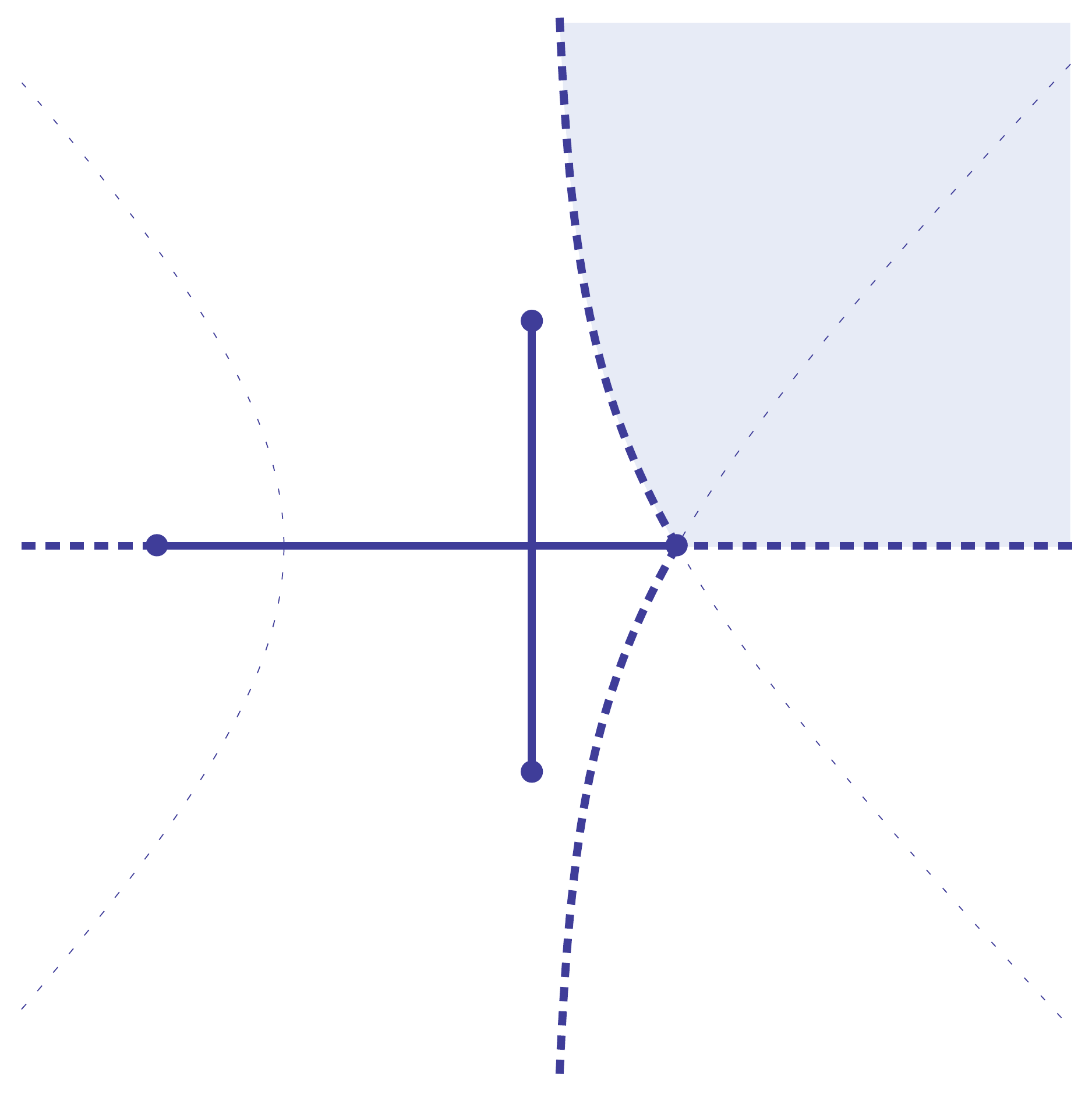}
      \put(65,75){$D_1$}
      \put(28,69){$\frac{c_2 + i|c_1|}{2\alpha}$}
      \put(28,28){$\frac{c_2 - i|c_1|}{2\alpha}$}
      \put(61,43){$K$}
      \put(6,43){$-3K$}
\end{overpic}
     \begin{figuretext}\label{cutsA.pdf}
       The domain $D_1$ for the triples given in (\ref{admissibleA}). The solid vertical segment is the branch cut $(-\frac{ic}{2\alpha}, \frac{i\bar{c}}{2\alpha})$; this cut collapses to a point when $c_1 = 0$. The solid horizontal segment is the branch cut connecting the triple zero $K = \sqrt{\frac{|\omega|}{12}}$ of $\Omega^2(k)$ with the simple zero at $-3K$.
      \end{figuretext}
     \end{center}
\end{figure}

For the family of triples given in (\ref{admissibleA}), the structure of the set $\im \Omega(k) = 0$ is shown in Figure \ref{cutsA.pdf}. In terms of the triple zero $K = \sqrt{\frac{|\omega|}{12}}$, we can write the family (\ref{admissibleA}) as follows:
\begin{align*}
\bigg\{\bigg(\alpha = \frac{8K^3}{c_2}, \omega = -12K^2, c & = \pm \sqrt{-c_2^2 + 12 K^4 + (\alpha^2 - 6 K^2)^2} + ic_2 
	\\
& \hspace{3cm}  \bigg| \;  0 < c_2 \leq 4K^2, \; K > 0\bigg\}.
\end{align*}
Note that $c_1 = 0$ iff $c_2 = 4K^2$.
Moreover, the branch cut $(-\frac{ic}{2\alpha}, \frac{i\bar{c}}{2\alpha})$ never intersects $\bar{D}_1$. Indeed, if $c_1 = 0$, then $-\frac{ic}{2\alpha} = \frac{i\bar{c}}{2\alpha} = K$, so the branch cut is absent. 
Moreover, since $\frac{c_2}{2\alpha} \leq K$, the cut always lies to the left of the triple zero. If $c_1 \neq 0$, then the value of $\im \Omega^2(k)$ is strictly positive at the top branch point,
$$\im \Omega^2\bigg(\frac{c_2 + i |c_1|}{2\alpha}\bigg) > 0,$$
showing that the cut does not intersect $\bar{D}_1$.

\begin{figure}
\begin{center}
\begin{overpic}[width=.45\textwidth]{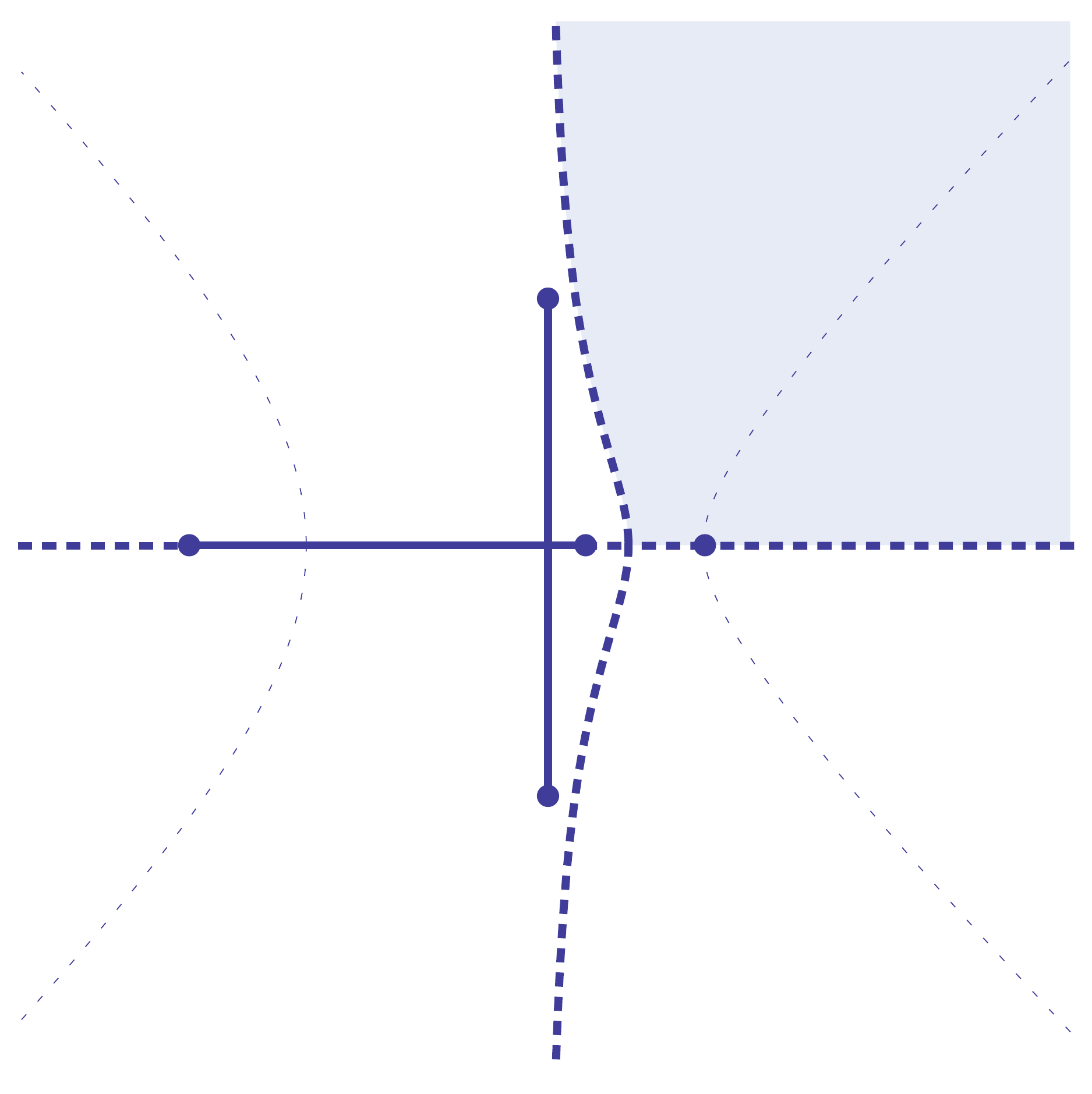}
      \put(65,75){$D_1$}
      \put(30,71){$\frac{c_2 + i|c_1|}{2\alpha}$}
      \put(30,26){$\frac{c_2 - i|c_1|}{2\alpha}$}
      \put(61,43){$K$}
      \put(6,43){$-3K$}
\end{overpic}
     \begin{figuretext}\label{cutsB.pdf}
       The domain $D_1$ for the triples given in (\ref{admissibleB}). 
      \end{figuretext}
     \end{center}
\end{figure}

\subsection{The family (\ref{admissibleB})}
For the family of triples given in (\ref{admissibleB}), the structure of the set $\im \Omega(k) = 0$ is shown in Figure \ref{cutsB.pdf}.
Since
$$c_1^2 = \frac{\left(c_2^2+8 K^4+2 K^2 \omega \right)^2 \left(\left(4 K^2+\omega
   \right)^2-4 c_2^2\right)}{4 c_2^4},$$
we find that $c_1 = 0$ iff $c_2 = -\frac{4K^2 + \omega}{2}$. 
Moreover, the branch cut $(-\frac{ic}{2\alpha}, \frac{i\bar{c}}{2\alpha})$ never intersects $\bar{D}_1$. Indeed, since $\frac{c_2}{2\alpha} \leq K$, the cut always lies to the left of $K$. Moreover, the value of $\im \Omega^2(k)$ is strictly positive at the top branch point,
$$\im \Omega^2\bigg(\frac{c_2 + i |c_1|}{2\alpha}\bigg) > 0,$$
showing that the cut does not intersect $\bar{D}_1$.

\subsection{The family (\ref{admissibleC})}
\begin{figure}
\begin{center}
\begin{overpic}[width=.45\textwidth]{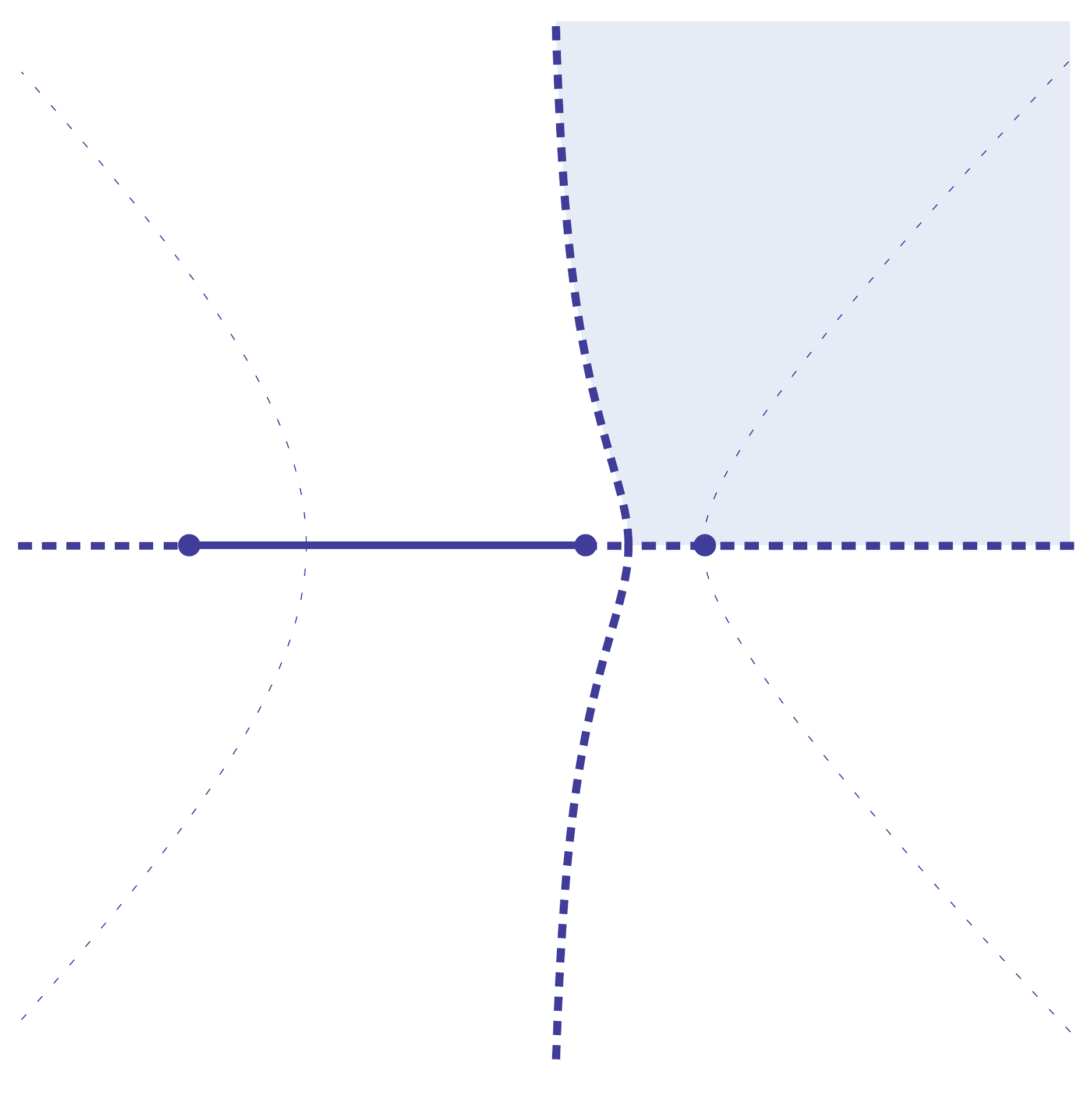}
      \put(65,80){$D_1$}
      \put(66,44){$K$}
   \put(4,44){$-K-\alpha$}
   \put(34,44){$-K+\alpha$}
\end{overpic}
     \begin{figuretext}\label{cutsC.pdf}
        The domain $D_1$ for the triples given in (\ref{admissibleC}). Since $c_1 = 0$, the branch cut $(-\frac{ic}{2\alpha}, \frac{i\bar{c}}{2\alpha})$ is absent. 
     \end{figuretext}
     \end{center}
\end{figure}

For the family of triples given in (\ref{admissibleC}), the structure of the set $\im \Omega(k) = 0$ is shown in Figure \ref{cutsC.pdf}.

In view of (\ref{case33}), we can write this family of triples as
\begin{align*}
&\bigg\{\bigg(\alpha  = \sqrt{\frac{|\omega|}{2} - 2K^2}, \omega, c = iK\sqrt{2 |\omega|-8K^2}\bigg)  
\; \bigg|  -12 K^2 < \omega < -4K^2, \; K > 0\bigg\}.
\end{align*} 
Since $c_1 = 0$, $-\frac{ic}{2\alpha} = \frac{i\bar{c}}{2\alpha} = K$.
The function $\Omega^2(k)$ has a double zero at $K$ and simple real zeros at
$-K \pm \alpha$. 

\subsection{The family (\ref{admissibleD})}
For the family of triples given in (\ref{admissibleD}), we distinguish three cases: $\omega > 0$, $\omega = 0$, and $-\alpha^2 \leq \omega < 0$. We consider each case in turn.

\subsubsection{$\omega > 0$}
\begin{figure}
\begin{center}
\begin{overpic}[width=.45\textwidth]{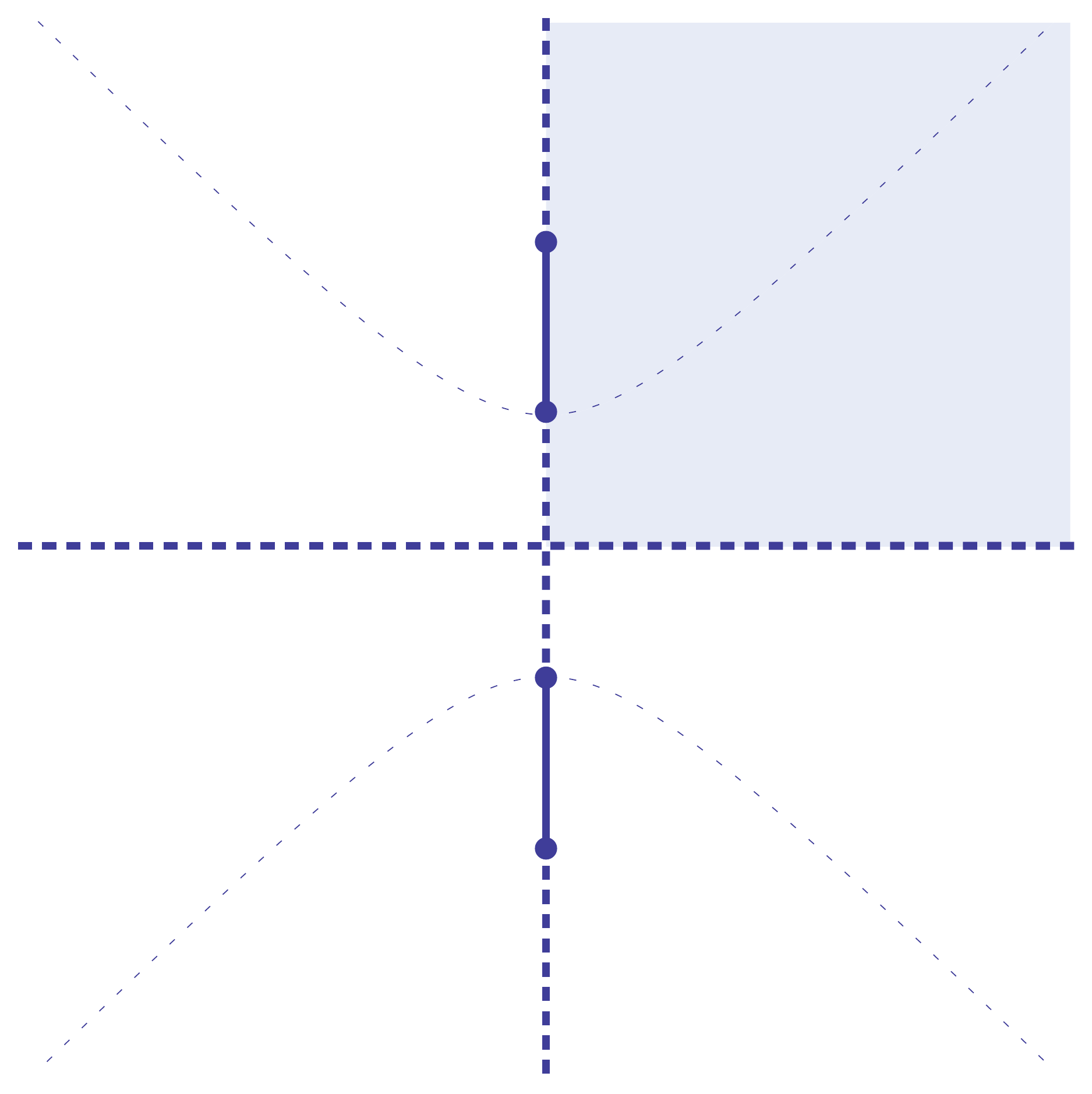}
      \put(75,65){$D_1$}
      \put(37,60){$\frac{i\sqrt{\omega}}{2}$}
      \put(32,37){$-\frac{i\sqrt{\omega}}{2}$}
      \put(52,78){$\frac{i\sqrt{\omega + \alpha^2}}{2}$}
      \put(52,20){$-\frac{i\sqrt{\omega + \alpha^2}}{2}$}
\end{overpic}
     \begin{figuretext}\label{cutsDa.pdf}
         The domain $D_1$ for the triples given in (\ref{admissibleD}) in the case of $\omega > 0$.    \end{figuretext}
     \end{center}
\end{figure}

The structure of the set $\im \Omega(k) = 0$ is shown in Figure \ref{cutsDa.pdf}.
There are four branch points, all purely imaginary, located at the double zeros of $\Omega^2(k)$ and at $(-\frac{ic}{2\alpha}, \frac{i\bar{c}}{2\alpha})$. More explicitly, the branch points are given by
$$\pm \frac{i\sqrt{\omega}}{2}, \pm \frac{i\sqrt{\omega + \alpha^2}}{2}.$$

\subsubsection{$\omega = 0$}
\begin{figure}
\begin{center}
\begin{overpic}[width=.45\textwidth]{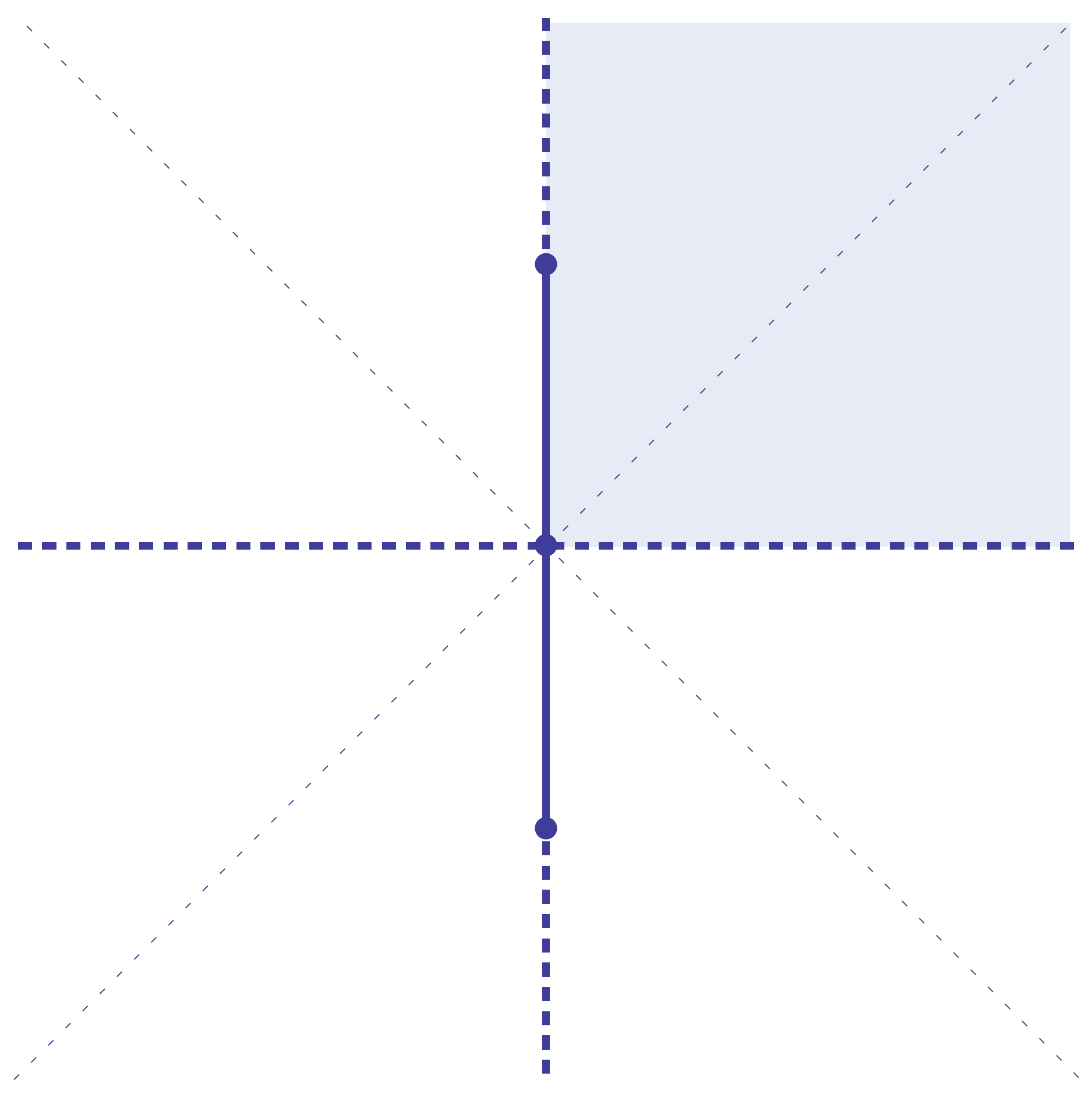}
      \put(75,65){$D_1$}
      \put(52,74){$\frac{i\alpha}{2}$}
      \put(52,23){$-\frac{i\alpha}{2}$}
\end{overpic}
     \begin{figuretext}\label{cutsDb.pdf}
         The domain $D_1$ for the triples given in (\ref{admissibleD}) in the case of $\omega = 0$.
     \end{figuretext}
     \end{center}
\end{figure}

The structure of the set $\im \Omega(k) = 0$ is shown in Figure \ref{cutsDb.pdf}.
There are three purely imaginary branch points at the fourth order zero of $\Omega^2(k)$ and at $(-\frac{ic}{2\alpha}, \frac{i\bar{c}}{2\alpha})$. The branch points are
$$0, \pm \frac{i \alpha}{2}.$$

\subsubsection{$-\alpha^2 \leq \omega < 0$}
\begin{figure}
\begin{center}
\begin{overpic}[width=.45\textwidth]{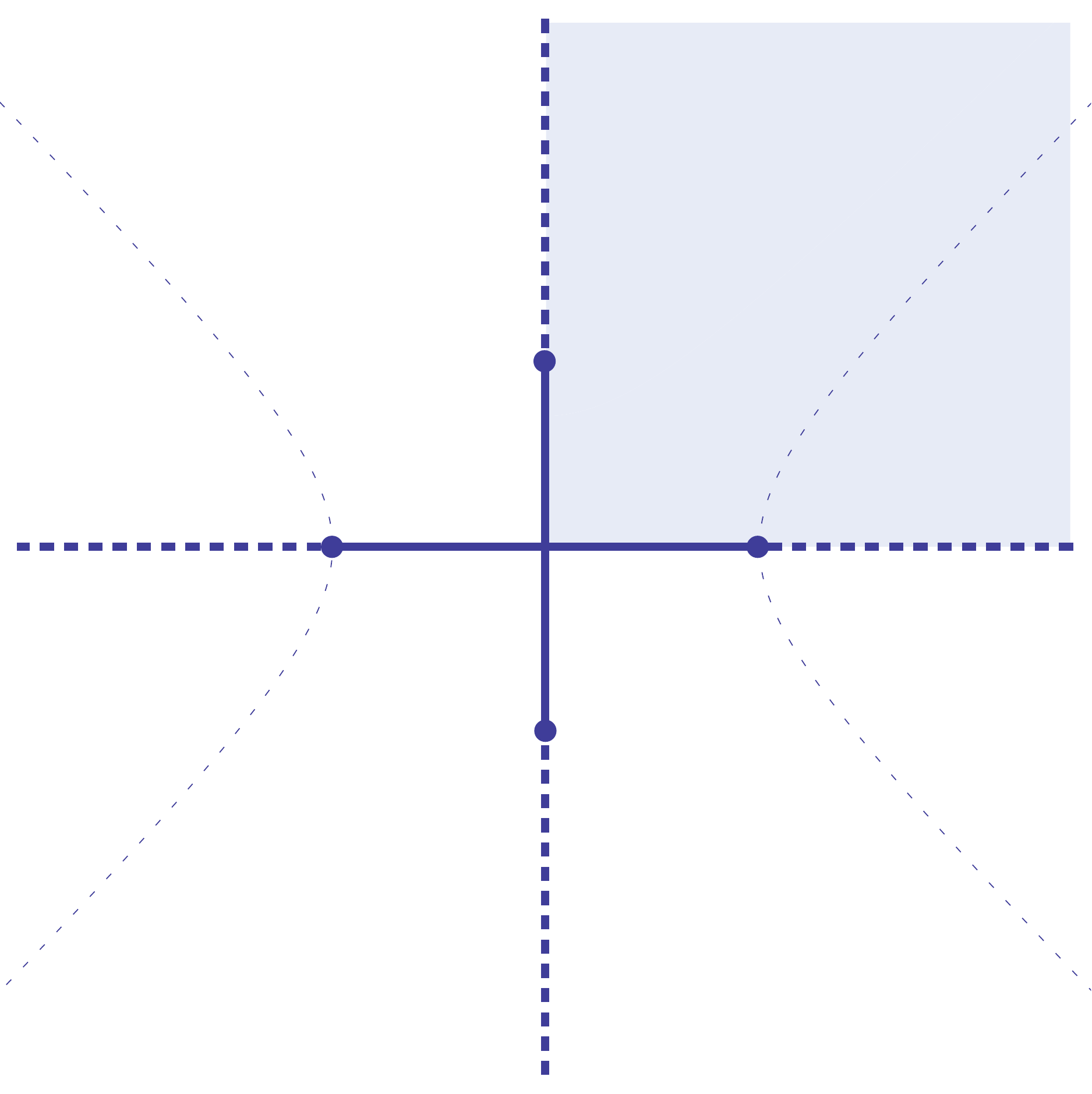}
     \put(75,80){$D_1$}
      \put(66,40){$\frac{\sqrt{|\omega|}}{2}$}
      \put(17,40){$-\frac{\sqrt{|\omega|}}{2}$}
      \put(52,67){$\frac{i\sqrt{\omega + \alpha^2}}{2}$}
      \put(51,29){$-\frac{i\sqrt{\omega + \alpha^2}}{2}$}
\end{overpic}
     \begin{figuretext}\label{cutsDc.pdf}
         The domain $D_1$ for the triples given in (\ref{admissibleD}) in the case of $-\alpha^2 \leq \omega < 0$.
      \end{figuretext}
     \end{center}
\end{figure}

The structure of the set $\im \Omega(k) = 0$ is shown in Figure \ref{cutsDc.pdf}.
There are four branch points at the double zeros of $\Omega^2(k)$ and at $(-\frac{ic}{2\alpha}, \frac{i\bar{c}}{2\alpha})$. The branch points are
$$\pm \frac{\sqrt{|\omega|}}{2}, \pm \frac{i\sqrt{\omega + \alpha^2}}{2}.$$

\subsection{The family (\ref{admissibleE})}
For the family of triples given in (\ref{admissibleE}), we distinguish two cases: $\omega = -3K^2$ and $-4K^2 < \omega < -3K^2$. 

\subsubsection{$\omega = -3K^2$}
\begin{figure}
\begin{center}
\begin{overpic}[width=.45\textwidth]{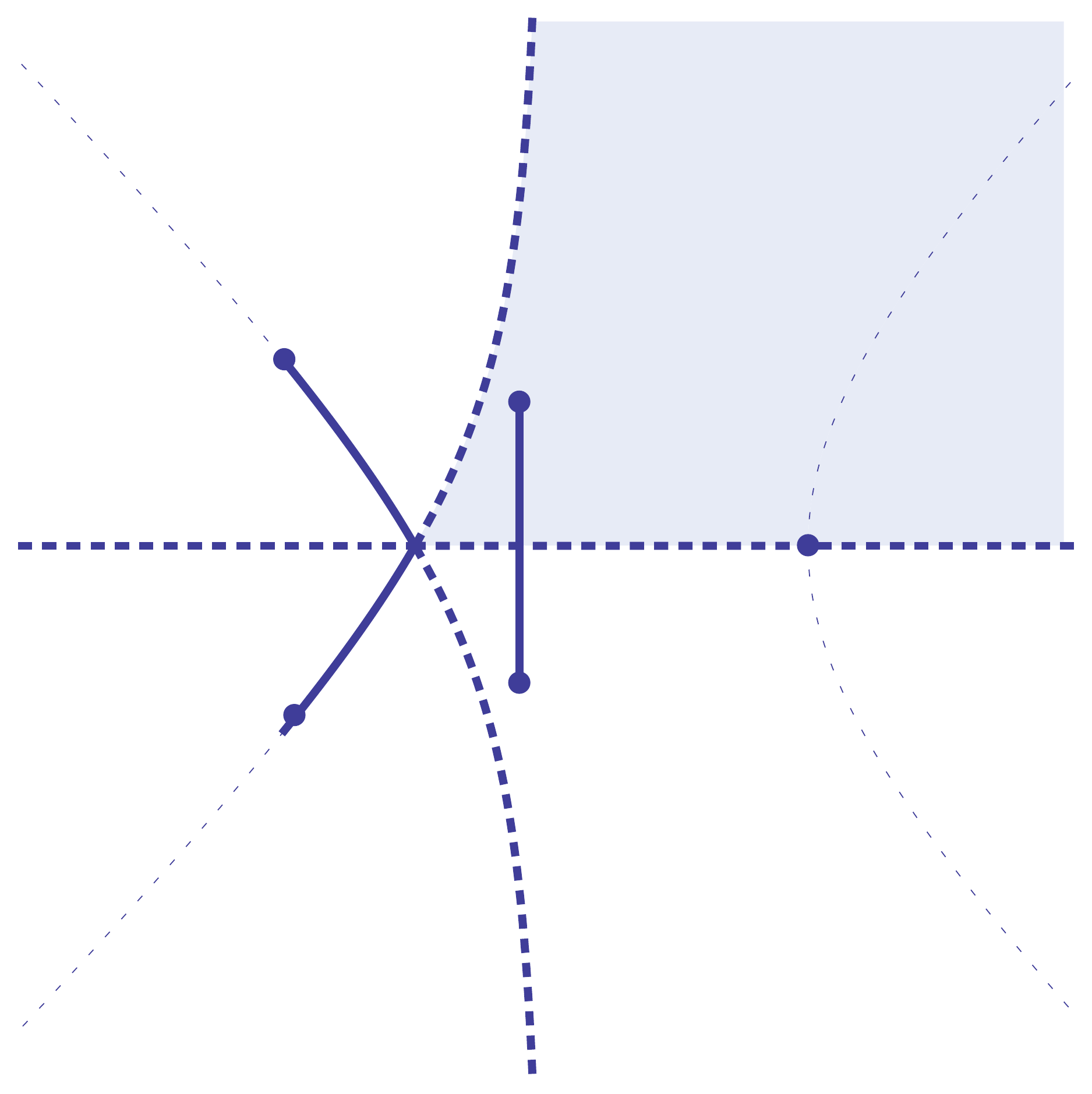}
     \put(65,75){$D_1$}
     \put(75,44){$K$}
      \put(-2,65){$-K + \frac{i K}{\sqrt{2}}$}
      \put(-2,32){$-K - \frac{i K}{\sqrt{2}}$}
      \put(50,63){$\frac{c_2 + i|c_1|}{2\alpha}$}
      \put(50,36){$\frac{c_2 - i|c_1|}{2\alpha}$}
\end{overpic}
     \begin{figuretext}\label{cutsEa.pdf}
         The domain $D_1$ for the triples given in (\ref{admissibleE}) in the case of $\omega = -3K^2$.
     \end{figuretext}
  \end{center}
\end{figure}

The structure of the set $\im \Omega(k) = 0$ is shown in Figure \ref{cutsEa.pdf}.
The function $\Omega^2(k)$ has a double zero at $K>0$ and two simple zeros at $(-1 \pm \frac{i}{\sqrt{2}})K$. 

Since
$$c_1^2 = -\frac{\left(4 c_2^2-K^4\right) \left(c_2^2+2 K^4\right)^2}{4c_2^4},$$
we see that $c_1 = 0$ iff $c_2 = -\frac{K^2}{2}$.
If $c_1 = 0$, the branch cut  $(-\frac{ic}{2\alpha}, \frac{i\bar{c}}{2\alpha})$ is absent. If $c_1 \neq 0$, the value of $\im \Omega^2(k)$ is strictly negative at the top branch point:
$$\im \Omega^2\bigg(\frac{c_2 + i |c_1|}{2\alpha}\bigg) < 0.$$
 Since $\frac{c_2}{2\alpha} < K$ the cut lies to the left of the double zero $K$ of $\Omega^2(k)$.
Hence, the cut $(-\frac{ic}{2\alpha}, \frac{i\bar{c}}{2\alpha})$ intersects $D_1$ if $c_1 \neq 0$.

\subsubsection{$-4K^2 < \omega < -3K^2$}
\begin{figure}
\begin{center}
\begin{overpic}[width=.45\textwidth]{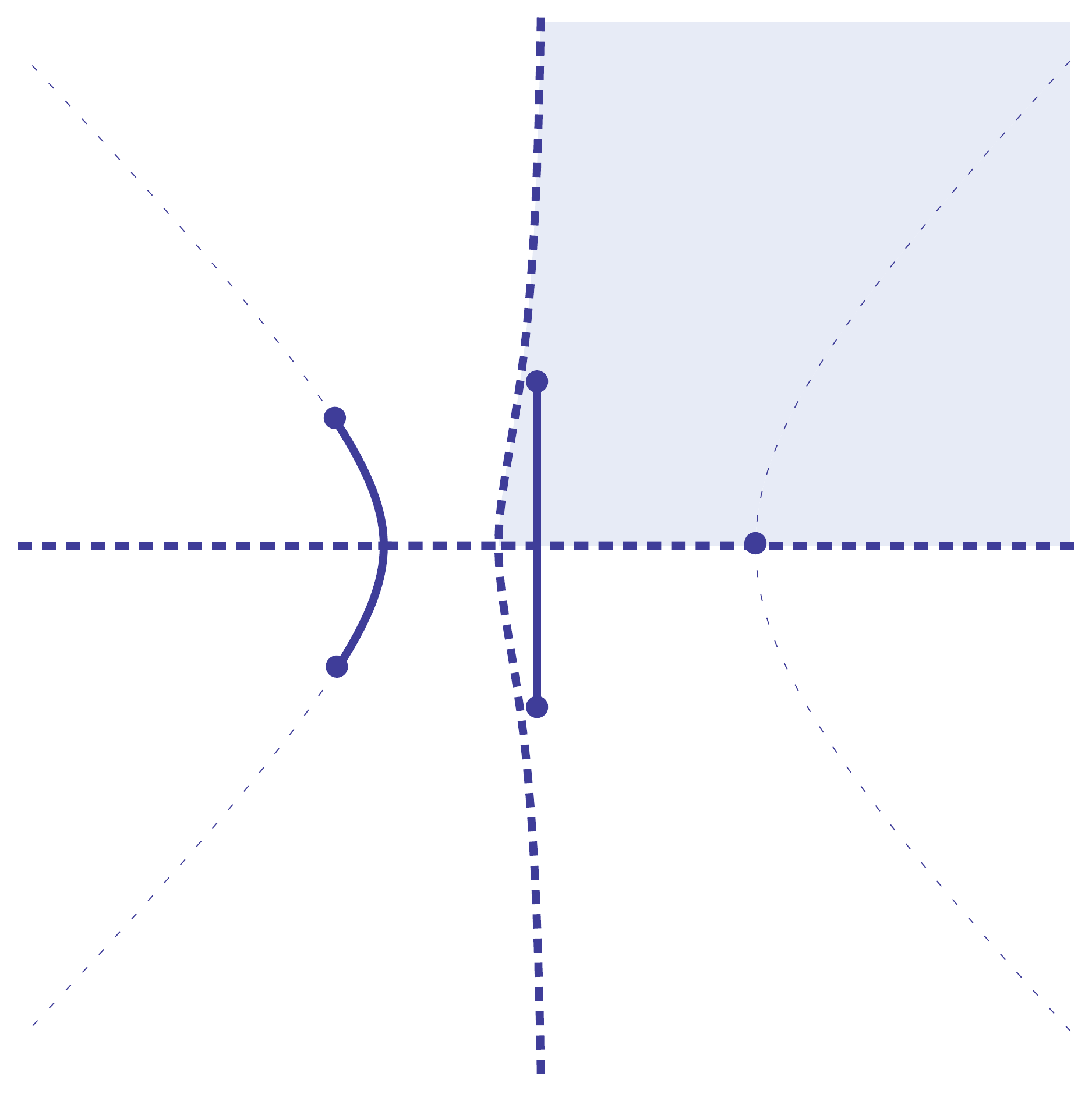}
     \put(65,75){$D_1$}
     \put(70,44){$K$}
%      \put(-2,65){$-K + \sqrt{-2K^2 - \frac{\omega}{2}}$}
 %     \put(-2,32){$-K - \sqrt{-2K^2 - \frac{\omega}{2}}$}
      \put(51,64){$\frac{c_2 + i|c_1|}{2\alpha}$}
      \put(51,33){$\frac{c_2 - i|c_1|}{2\alpha}$}
\end{overpic}
     \begin{figuretext}\label{cutsEb.pdf}
         The domain $D_1$ for the triples given in (\ref{admissibleE}) in the case of $-4K^2 < \omega < -3K^2$. 
     \end{figuretext}
  \end{center}
\end{figure}

The structure of the set $\im \Omega(k) = 0$ is shown in Figure \ref{cutsEb.pdf}.
The function $\Omega^2(k)$ has a double zero at $K>0$ and two simple zeros at $-K \pm \sqrt{-2K^2 - \frac{\omega}{2}}$. We have $c_1 = 0$ iff $c_2 = -\frac{4K^2 + \omega}{2}$.

The branch cut $(-\frac{ic}{2\alpha}, \frac{i\bar{c}}{2\alpha})$ is a vertical segment intersecting the real axis at $\frac{c_2}{2\alpha}$. If $c_1 = 0$, this cut is absent. Assume $c_1 \neq 0$. Since $\frac{c_2}{2\alpha} < 0$ the cut lies in the left half-plane. 
Let $K_1 < K_2 < 0 <  K$ denote the three zeros of the polynomial $P(k) = 4k^3 + \omega k + \alpha c_2 = 0$, i.e., $K_1, K_2, K$ are the three intersection points of $\Gamma$ with the real axis. Since $P(\frac{c_2}{2\alpha}) < 0$, we have $\frac{c_2}{2\alpha}<K_1$ or $K_2 < \frac{c_2}{2\alpha} < 0$. Since $P'( \frac{c_2}{2\alpha}) < 0$, we have in fact $K_2 < \frac{c_2}{2\alpha} < 0$. 
Moreover, the value of $\im \Omega^2(k)$ is strictly negative at the top branch point:
$$\im \Omega^2\bigg(\frac{c_2 + i |c_1|}{2\alpha}\bigg) < 0.$$
This shows that the branch cut $(-\frac{ic}{2\alpha}, \frac{i\bar{c}}{2\alpha})$ intersects $D_1$ and has the qualitative form shown in Figure \ref{cutsEb.pdf} for $c_1 \neq 0$.

%For the defocusing NLS on the half-line  $a(k)$  cannot vanish in $\im k \geq 0$, see \cite{BK2000} p. 1820. (For $k \in \R$ this is clear from the determinant relation.)

\bigskip
\noindent
{\bf Acknowledgement} {\it The author acknowledges support from the EPSRC, UK.}

\bibliographystyle{plain}
\bibliography{is}

\end{document}